\documentclass[smallextended]{svjour3}

\smartqed  

\usepackage{amssymb}
\usepackage{empheq}
\usepackage[mathscr]{euscript}
\usepackage{stmaryrd}
\usepackage{mathabx}
\usepackage{dsfont}
\usepackage{graphicx}
\usepackage{multirow}
\usepackage{algorithm2e}
\usepackage[table]{xcolor}
\usepackage{color}
\usepackage{tikz}
\usetikzlibrary{3d}
\usetikzlibrary{decorations.pathreplacing,angles,quotes}
\usepackage{pgfplots}
\usepackage{algorithm2e}

\usepackage[colorlinks=true,breaklinks=true,bookmarks=true,urlcolor=blue,
     citecolor=blue,linkcolor=blue,bookmarksopen=false,draft=false]{hyperref}

\newcommand{\calA}{\mathcal{A}}

\newcommand{\calL}{\mathcal{L}}
\newcommand{\calP}{\mathcal{P}}
\newcommand{\calT}{\mathcal{T}}

\newcommand{\scrO}{\mathscr{O}}

\newcommand{\bbN}{\mathbb{N}}
\newcommand{\bbR}{\mathbb{R}}
\newcommand{\bbZ}{\mathbb{Z}}

\newcommand{\aff}{\operatorname{aff}}

\newcommand{\Conv}{\operatorname{Conv}}
\newcommand{\Em}{\operatorname{Em}}
\newcommand{\ext}{\operatorname{ext}}

\newcommand{\Proj}{\operatorname{Proj}}
\newcommand{\relint}{\operatorname{relint}}
\newcommand{\rev}{\operatorname{rev}}
\newcommand{\Span}{\operatorname{span}}
\newcommand{\Slice}{\operatorname{Slice}}
\newcommand{\SOS}{\operatorname{SOS2}}

\newcommand{\supp}{\operatorname{supp}}

\newcommand\defeq{\mathrel{\overset{\makebox[0pt]{\mbox{\normalfont\tiny\sffamily def}}}{=}}}
\newcommand{\bra}[1]{\left(#1\right)}

\renewcommand{\int}{\operatorname{int}}

\newcommand{\polylog}{\operatorname{polylog}}

\newcommand{\HMC}{\ensuremath{H^{\operatorname{mc}}}}
\newcommand{\HEX}{\ensuremath{H^{\operatorname{ex}}}}
\newcommand{\HRB}{\ensuremath{H^{\operatorname{br}}}}
\newcommand{\HZZ}{\ensuremath{H^{\operatorname{zz}}}}
\newcommand{\TSOS}{\ensuremath{\calT^{\operatorname{SOS2}}}}
\newcommand{\TANN}{\ensuremath{\calT^{\operatorname{ann}}}}

\begin{document}

\title{A mixed-integer branching approach for very small formulations of disjunctive constraints}

\titlerunning{Mixed-integer branching formulations}

\author{Joey Huchette \and Juan Pablo Vielma}

\institute{J. Huchette \at \email{huchette@mit.edu} \and J. P. Vielma \at \email{jvielma@mit.edu}}

\maketitle

\begin{abstract}

   An important problem in optimization is the construction of mixed-integer programming (MIP) formulations of disjunctive constraints that are both strong and small. Motivated by lower bounds on the number of integer variables that are required by traditional MIP formulations, we present a more general mixed-integer branching formulation framework. Our approach maintains favorable algorithmic properties of traditional MIP formulations: in particular, amenability to branch-and-bound and branch-and-cut algorithms. Our main technical result gives an explicit linear inequality description for both traditional MIP and mixed-integer branching formulations for a wide range of disjunctive constraints. The formulations obtained from this description have linear programming relaxations that are as strong as possible and generalize some of the most computationally effective formulations for piecewise linear functions and other disjunctive constraints.  We use this result to produce a strong mixed-integer branching formulation for any disjunctive constraint that uses only two integer variables and a linear number of extra constraints. We sharpen this result for univariate piecewise linear functions and annulus constraints arising in power systems and robotics, producing strong mixed-integer branching formulations that use only two integer variables and a constant ($\leq 6$) number of general inequality constraints. Along the way, we produce two strong logarithmic-sized traditional MIP formulations for the annulus constraint using our main technical result, illustrating its broader utility in the traditional MIP setting.

\end{abstract}



\section{Introduction} \label{sec:introduction}

Consider a disjunctive set $S = \bigcup_{i=1}^d P^i$, where each $P^i \subset \bbR^n$ is a rational bounded polyhedra. Disjunctive constraints of the form $x \in S$ abound in optimization: they are useful, for example, to model nonlinearities~\cite{Geisler:2012,Misener:2012} or discrete logic imposed by complex processes~\cite{Bergamini:2005,Trespalacios:2014b}. Therefore, we would like a way to represent these constraints in such a way that we can efficiently optimize over them. Additionally, we would like to do this in a composable way, as disjunctive constraints frequently arise as substructures in large, complex optimization problems.

Mixed-integer programming (MIP) offers one such solution. A MIP formulation for $S$ is given by a linear programming (LP) relaxation
\[
    R = \left\{(x,y,z) \in \bbR^{n+m+r} : Ax + By + Cz \leq d \right\}
\]
such that, when integrality is imposed on the integer (control) variables $z$, the set projects down onto $S$. MIP formulations are useful because there are sophisticated algorithms--and corresponding high-quality software implementations--that can optimize over these representations efficiently in practice~\cite{Bixby:2007,Junger:2010}. Furthermore, combining MIP formulations for different substructures is trivial, and so this technology can be marshalled for very complex and large-scale optimization problems.

Indeed, it is very often the case that MIP formulations outperform branch-and-bound methods that work directly on the disjunctions (e.g. \cite{Beaumont:1990,Land:1960}), despite the fact that they require additional integer variables and constraints. This can largely be attributed to the immense advances of MIP solvers over the past decades, as their algorithms are now much more complex than a traditional branch-and-bound approach. In particular, the development of a sophisticated theory on cutting planes plays a crucial role here~\cite{Bixby:2007,Junger:2010}. These techniques can easily combine information from multiple disjunctions and other constraints in the optimization problem to provide tighter relaxations and shorten computation time. Achieving this without the integer control variables can require significant theoretical developments even for very specific structures~\cite{Farias-Jr.:2013,Keha:2006}. The downside of the MIP approach, however, is that it requires additional integer control variables and constraints, which leads to larger (and therefore slower) LP relaxations.

One goal of this work is to reduce the number of control variables needed to model disjunctions. A folklore result holds that any MIP formulation must use at least $\lceil \log_2(d) \rceil$ control variables to represent the union of $d$ sets. A recent generalization of this result~\cite[Lemma 2]{Lubin:2017a} shows that this bound holds even if we allow general integer control variables or a convex nonlinear relaxation. Therefore, if we hope to further reduce the number of control variables, we will need to go beyond traditional MIP formulations. 

Fortunately, there has recently been growing interest in studying the expressive power and computational properties of generalizations of traditional MIP formulations \cite{Bader:2015,Bonami:2017,Hildebrand:2017}. In particular, our work builds off of the ideas of Bonami et al.~\cite{Bonami:2017} for handling ``holes'' in integer sets. As a simple example, consider the disjunctive constraint $x \in S = \{1,2,4,5\}$. This constraint is nearly equivalent to a standard integrality constraint $x \in [1,5] \cap \bbZ$, but with a hole in the domain at $3$. A traditional MIP formulation for this might introduce a binary variable $z$ and impose the constraints
\begin{equation} \label{eqn:simple-MIP-formulation}
    1z + 4(1-z) \leq x \leq 2z + 5(1-z), \quad\quad (x,z) \in \bbZ \times \{0,1\}.
\end{equation}

Bonami et al. handle these holes directly in the original $x$ space, using \emph{wide split disjunctions}. A standard branch-and-bound algorithm~\cite{Land:1960} will perform \emph{variable branching} on fractional solutions: given the point $\hat{x}=2.5$, it rounds $\hat{x}$ and imposes the valid disjunction $x \leq 2 \vee x \geq 3$ to separate this point. However, this leaves us no way to separate the hole at $\hat{x}=3$, which is integer but not feasible for the original constraint. A natural way around this is to impose a wide split disjunction of the form $x \leq 2 \vee x \geq 4$ to separate $\hat{x}$ from $S$. This is a straightforward change to the branch-and-bound algorithm in a way that does not require any additional control variables. The crucial observation of Bonami et al. is that wide split disjunctions also readily admit standard cutting plane techniques, such as the intersection cut~\cite{Andersen:2005}. By combining the slightly modified branch-and-bound algorithm with cutting planes, Bonami et al. observe a considerable computational speed-up compared to a ``full'' formulation like \eqref{eqn:simple-MIP-formulation} when optimizing over integer sets with holes.

Our work extends this general idea to the mixed-integer setting with the aim of constructing very small formulations for disjunctive constraints. We will see that if we allow holes in our mixed-integer formulations, we can drastically reduce the number of control variables and constraints we need to build formulations. We can optimize over these representations using variable branching or wide split disjunctions, meaning that the same cutting plane machinery (and, hopefully, computational performance) applied by Bonami et al. are applicable in our case as well. In certain degenerate cases we will need to deploy two-term non-parallel disjunctions, for which cut generation techniques have also been developed in the literature~\cite{Andersen:2005,Balas:1998,Bonami:2013,Kis:2014}. 

More concretely, our contributions are as follows.

\begin{itemize}
    \item \textbf{An explicit geometric construction for strong formulations of disjunctive constraints.} This gives us a practical way to construct both traditional and generalized MIP formulations for the broad class of \emph{combinatorial disjunctive constraints}~\cite{Huchette:2016a}. For the traditional case, the resulting formulations are \emph{integral} or \emph{ideal} (i.e. the linear programming relaxation of the formulations have extreme points that satisfy the integrality constraints on the control variables), and for the generalized case they satisfy a natural generalization of this property. The construction also gives an upper bound on the number of general inequality constraints needed to construct an ideal formulation (traditional and/or generalized) for a given constraint. 
    \item \textbf{A framework for generalized MIP formulations, and branching rules to optimize over them.} We present the family of \emph{mixed-integer branching formulations} as a mixed-integer generalization of the ``integer programming with holes'' approach of Bonami et al~\cite{Bonami:2017}. These generalized MIP formulations sidestep the logarithmic lower bound on the number of integer variables and have practical branching rules that can be used to generate cutting planes and implement branch-and-bound algorithms. Finally, explicit inequality descriptions for these formulations can be easily obtained from our main geometric construction result.
    \item \textbf{Very small formulations for disjunctive constraints.} We show that for any combinatorial disjunctive constraint there exists an ideal mixed-integer branching formulation with only two control variables and at most a linear number of constraints. For the SOS2 constraint, we improve this with an ideal mixed-integer branching formulation using two control variables and only four general inequality constraints. We also study a relaxation for the annulus, for which we construct an ideal mixed-integer branching formulation with two control variables and only six general inequality constraints. Finally, we apply our main result to also produce two new logarithmic-sized traditional MIP formulations for this relaxed annulus constraint.
\end{itemize}

\section{Preliminaries}

We will use the following generalization of the traditional MIP formulation introduced in Section~\ref{sec:introduction}.

\begin{definition} \label{def:valid-formulation}
     Take some set $S \subseteq \bbR^n$, along with a rational polyhedron $R = \left\{(x,z) \in \bbR^{n+r} : Ax + Cz \leq d\right\}$.
    \begin{itemize}
        \item  We say that $F = \{\bra{x,z}\in R : z\in H\}$ is a \emph{formulation} of $x\in S$ (or just $S$) with respect to the set $H \subset \bbR^r$ if $\Proj_x\bra{F}=S$.
            \item We refer to $R$ as the \emph{linear programming (LP) relaxation} of the formulation.

            \item We call $x$ the \emph{original variables} and $z$ the \emph{control variables}.

    \item We say a linear inequality  defining $R$ is a \emph{variable bound} if it has only one non-zero coefficient, and a \emph{general inequality constraint} otherwise.
    \end{itemize}
\end{definition}

Note that we are omitting auxiliary continuous variables in our formulation description. In theory, auxiliary variables could drastically reduce the number of constraints needed to describe the LP relaxation $R$. However, the cases we study in this work will admit very small formulations without auxiliary variables, so we omit them from our discussion for clarity. We also differentiate the two types of constraints describing our relaxation as a simplex-based algorithm can typically impose variable bounds with little or no computational overhead. Finally, we use the ``control variables'' terminology of Jeroslow~\cite{Jeroslow:1987} to emphasize that these variables will not necessarily be allowed to take arbitrary integer values.

We can easily recover traditional MIP formulations as a special case of our more general definition.
\begin{definition}
    We say that $H \subseteq \bbZ^r$ is \emph{hole-free} if $\Conv(H) \cap \bbZ^r = H$.
\end{definition}
If $H$ is hole-free, we can replace the set constraint $z \in H$ in our formulation with $z \in \Conv(H) \cap \bbZ^r$. In the case where $\Conv(H)$ is a polyhedron (which will be the case for the remainder), we recover a traditional (linear) MIP formulation for $S$. Of particular note are \emph{binary MIP formulations}, which correspond to the case where $H \subseteq \{0,1\}^r$.

The usual notions of formulation strength~\cite{Vielma:2015} also carry over directly to our more general setting.

\begin{definition}
    A formulation of $S \subset \bbR^n$ with respect to $H$ is \emph{ideal} if the extreme points of its LP relaxation $R$ satisfy $\ext(R) \subseteq \bbR^{n} \times H$.
\end{definition}

\subsection{The embedding approach}

We will construct formulations for disjunctive sets $S = \bigcup_{i=1}^d P^i$ through what is known as the embedding approach~\cite{Vielma:2015a}. We assign each \emph{alternative} $P^i$ a unique \emph{code} $h^i \in \bbR^r$. We call such a collection of distinct vectors $H = (h^i)_{i=1}^d$ an \emph{encoding}. Given $\calP = (P^i)_{i=1}^d$ and $H$, we construct the \emph{embedding} of $S$ in a higher-dimensional space as
\[
    \Em(\calP,H) \defeq \bigcup_{i=1}^d (P^i \times \{h^i\}).
\]
This object is useful as projecting out the control variables gives us the disjunctive set: $\Proj_x(\Em(\calP,H)) = S$. In particular, if the encoding satisfies a natural geometric condition, then $Q(\calP,H) \defeq \Conv(\Em(\calP,H))$ immediately gives us the LP relaxation for an ideal formulation of $S$ with respect to $H$. 

\begin{definition}
    A set $H \subset \bbR^r$ is in \emph{convex position} if $\ext(\Conv(H)) = H$.
\end{definition}

\begin{proposition}
    Take $\calP = (P^i)_{i=1}^d$. If $H=(h^i)_{i=1}^d$ is an encoding in convex position, then $\{(x,z) \in Q(\calP,H) : z \in H\}$ is an ideal formulation for $\bigcup_{i=1}^d P^i$ with respect to $H$.
\end{proposition}
To use an embedding formulation like this in practice we will need (1) an explicit outer (inequality) description of $Q(\calP,H)$, and (2) a way to iteratively impose the set constraint $z\in H$. Theorem~\ref{thm:general-cdc-characterization} will give use a way to meet this first requirement in both the traditional and generalized MIP setting. For the second requirement, variable branching suffices in the traditional MIP setting when $H \subseteq \bbZ^r$ is hole-free, as $z\in \Conv(H) \cap \bbZ^r = H$. We will develop analogous branching schemes for our generalized MIP formulations in Section~\ref{section3}.

\subsection{Combinatorial disjunctive constraints}

In this work, we will be primarily interested in constructing formulations for \emph{combinatorial disjunctive constraints}~\cite{Huchette:2016a}, where each alternative is some face on the unit simplex. Notationally, take:
\begin{itemize}
    \item $\llbracket d \rrbracket \defeq \{1,2,\ldots,d-1,d\}$ and $\llbracket j,k \rrbracket \defeq \{j,j+1,\ldots,k-1,k\}$;
    \item $[d]^2 \defeq \{\{i,j\} \in \llbracket d \rrbracket^2 : i < j\}$ (where notationally $\{i,j\} \in [d]^2$ implies $i < j$);
    \item The unit simplex as $\Delta^n \defeq \{\lambda \in \bbR^n_+ : \sum_{v=1}^n \lambda_v = 1\}$;
    \item The support of an element $\lambda \in \Delta^n$ as $\supp(\lambda) \defeq \{v \in \llbracket n \rrbracket : \lambda_v \neq 0\}$; and
    \item The face of $\Delta^n$ induced by $T \subseteq \llbracket n \rrbracket$ as $P(T) \defeq \{\lambda \in \Delta^n : \supp(\lambda) \subseteq T\}$.
\end{itemize}

\begin{definition}
    A \emph{combinatorial disjunctive constraint} is a constraint of the form $\lambda \in \bigcup_{i=1}^d P(T^i)$, given by the family of distinct nonempty sets $\calT = (T^i \subseteq \llbracket n \rrbracket)_{i=1}^d$. We will denote the corresponding set of alternatives as $\calP(\calT) \defeq (P(T^i))_{i=1}^d$.
\end{definition}

Throughout we will assume that $\bigcup_{i=1}^d T^i = \llbracket n \rrbracket$; or, equivalently, that $\Conv(\bigcup_{i=1}^d P(T^i)) = \Delta^n$. This is without loss of generality (w.l.o.g.), as otherwise we could simply drop any missing component from the constraint.

Due to the classical Minkowski-Weyl Theorem (e.g. \cite[Corollary 3.14]{Conforti:2014}), we can formulate \emph{any} disjunctive constraint as a combinatorial disjunctive constraint, provided we have a description for each alternative in terms of its extreme points. In particular, if we label the extreme points as $\bigcup_{i=1}^d \ext(P^i) = \{v^j\}_{j=1}^n$, then we can write the disjunctive set as
\[
    \bigcup_{i=1}^d P^i = \left\{ \sum_{j=1}^n \lambda_j v^j : \lambda \in \bigcup_{i=1}^d P(T^i) \right\},
\]
where each $T^i = \{j \in \llbracket n \rrbracket : v^j \in \ext(P^i)\}$ corresponds to the indices of the extreme points of $P^i$ in our ordering. Combinatorial disjunctive constraints are a particularly natural way to formulate a number of disjunctive constraints of interest, including piecewise linear functions, non-convex set inclusion or collision avoidance constraints, and relaxations for multilinear functions~\cite{Huchette:2016a}.

Unfortunately, combinatorial disjunctive constraints are a class of constraints for which the folklore lower bound on the number of integer control variables holds under a simple non-redundancy property. We can formalize this through the following simple proposition that we prove in Section~\ref{app:prove-cdc-properties}.

\begin{proposition}\label{easycombinatorial}
    Take $\calT = (T^i \subseteq \llbracket n \rrbracket)_{i=1}^d$ as a representation for a combinatorial disjunctive constraint, where $T^i \not\subseteq T^j$ and $T^j \not\subseteq T^j$ for each $\{i,j\} \in [d]^2$. Consider some encoding $H = (h^i)_{i=1}^d \subset \bbR^r$. Then there exists a formulation for $\bigcup_{i=1}^d P(T^i)$ with respect to $H$ if and only if $H$ is in convex position. In particular, if such a formulation exists and $H$ is hole-free, then $r \geq \lceil \log_2(d) \rceil$ necessarily.
\end{proposition}

Proposition~\ref{easycombinatorial} tells us that any traditional MIP formulation for a combinatorial disjunctive constraint requires at least a logarithmic number of control variables. Furthermore, Huchette and Vielma~\cite{Huchette:2017,Vielma:2015a} give a lower bound of $2\lceil \log_2(d) \rceil$ on the number of general inequality constraints for ideal (non-extended) formulations of the SOS2 constraint~\cite{Beale:1970}, where $\TSOS_d \defeq (\{i,i+1\})_{i=1}^d$. This result follows from an explicit outer description for $Q(\calP(\TSOS_d),H)$ for the case where $H$ is any hole-free encoding that is in convex position. In other words, the result characterizes all non-extended ideal MIP formulations of the SOS2 constraint. This result can be used to produce the family of ``logarithmic'' MIP formulations for the SOS2 constraint which have proven computational efficacy~\cite{Huchette:2017,Vielma:2015,Vielma:2010,Vielma:2009a}.

\subsection{Summary of main contributions}
Motivated by the success of the logarithmic traditional MIP formulations for the SOS2 constraint, we extend the geometric characterization of Huchette and Vielma~\cite{Huchette:2017,Vielma:2015a} to any combinatorial disjunctive constraint. More precisely, we give an explicit linear inequality description of $Q(\calP(\calT),H)$ for any combinatorial disjunctive constraint given by the family $\calT$, paired with any encoding $H$ that is in convex position. This gives a practical way to build ideal formulations, particularly for low-dimensional encodings. It also gives an upper bound on the number of general inequality constraints needed to construct any ideal formulation for a given combinatorial disjunctive constraint. The statement of our main technical result is as follows.

\begin{theorem} \label{thm:general-cdc-characterization}
    Take $\calT = (T^i \subseteq \llbracket n \rrbracket)_{i=1}^d$ and $H = (h^i)_{i=1}^d \subset \bbR^r$ as an encoding in convex position. Let $D = \{\{i,j\} \in [d]^2 : T^i \cap T^j \neq \emptyset\}$, and presume that $D$ is connected in the sense that the associated graph $G=(\llbracket d \rrbracket,D)$ is connected. Take $C = \{c^{i,j} \defeq h^j - h^i\}_{\{i,j\} \in D}$, and $\calL = \Span(C)$. Define $M(b;\calL) \defeq \{y \in \calL : b \cdot y = 0\}$ to be the hyperplane in the linear space $\calL$ induced by the direction $b \neq {\bf 0}^r$. If $\{b^k\}_{k=1}^\Gamma\subset \bbR^r \backslash \{{\bf 0}^r\}$ is such that $\{M(b^k;\calL)\}_{k=1}^\Gamma$ is the set of linear hyperplanes spanned by $C$ in $\calL$, then $(\lambda,z) \in Q(\calP(\calT),H)$ if and only if
    \begin{subequations} \label{eqn:general-V-formulation}
    \begin{gather}
        \sum_{v=1}^n \min_{s : v \in T^s}\{b^k \cdot h^s\} \lambda_v \leq b^k \cdot z \leq \sum_{v=1}^n \max_{s : v \in T^s}\{b^k \cdot h^s\} \lambda_v \quad \forall k \in \llbracket \Gamma \rrbracket \label{eqn:general-V-formulation-1} \\
        (\lambda,z) \in \Delta^{n} \times \aff(H). \label{eqn:general-V-formulation-2}
    \end{gather}
    \end{subequations}
\end{theorem}

We defer the proof to Section~\ref{sec:prove-cdc-characterization}, and instead concentrate on its implication for non-traditional formulations, which can be summarized as follows.
\begin{enumerate}
    \item \textbf{[Proposition~\ref{prop:general-cdc-moment-curve}]} For any combinatorial disjunctive constraint on $n$ components with $d$ alternatives, we can produce an ideal mixed-integer branching formulation with two control variables, $\scrO(d)$ general linear inequality constraints, $\scrO(n)$ variable bounds, and one equation.

    \item \textbf{[Proposition~\ref{prop:sos2-constant}]} For the SOS2 constraint on $n=d+1$ components, we can produce an ideal mixed-integer branching formulation with two control variables, four general linear inequality constraints, $\scrO(n)$ variable bounds, and one equation.

    \item \textbf{[Propositions~\ref{prop:log-annulus} and \ref{prop:zig-zag-annulus}]} For a relaxation of the annulus as a partition of $d$ quadrilaterals, we can produce two ideal traditional MIP formulations with $\lceil\log_2(d)\rceil$ control variables, $\scrO(\polylog(d))$ general linear inequality constraints, $\scrO(d)$ variable bounds, and one equation.

    \item \textbf{[Proposition~\ref{prop:exotic-annulus}]} For a relaxation of the annulus as a partition of $d$ quadrilaterals, we can produce an ideal mixed-integer branching formulation with two control variables, six general linear inequality constraints, $\scrO(d)$ variable bounds, and one equation.
\end{enumerate}

In other words, our new formulation approach allows us to construct ideal formulations for any combinatorial disjunctive constraint with very few control variables and at most a linear number of general linear inequality constraints.  Furthermore, by taking advantage of structure, we can further reduce this to a constant number of general inequalities for the SOS2 constraint and a tight relaxation for the annulus.

However, a formulation with two control variables implies a two-dimensional encoding, which cannot be hole-free if $d > 4$. Hence, the resulting formulation is not a traditional MIP formulation. In the following section, we present a way to optimize over such representations in a branch-and-bound setting, using branching schemes customized for a particular encoding. The branching schemes we present will use combinations of variable branching, wide axis-aligned split disjunctions (or just \emph{wide variable branching}) a la Bonami et al.~\cite{Bonami:2017}, and general two-term disjunctions. As a result, both standard and state-of-the-art cutting plane technology can be deployed to strengthen the relaxations of our formulations.

\section{Branching schemes and mixed-integer branching formulations} \label{section3}

We reiterate that traditional MIP formulations are useful because there exist algorithms--and high-quality implementations of those algorithms--that are able to optimize over these representations efficiently in practice. Roughly, these implementations work by applying the branch-and-bound method~\cite{Land:1960}, coupled with the judicious application of cutting planes to strengthen the LP relaxation. In this section, we formalize how our generalized notion of a MIP formulation fits into the branch-and-bound framework.

\subsection{Branching schemes}
We start by formally defining what we mean by a branching scheme.

\begin{definition} \label{def:branching-scheme}
    A \emph{branching scheme} is a procedure that, given
    \begin{itemize}
        \item a polyhedron $Q \subset \bbR^{r}$,
        \item an encoding $H \subset \bbR^r$, and
        \item a point $\hat{z} \in Q$,
    \end{itemize}
    either verifies that $\hat{z} \in H$, or outputs two polyhedra $Q^1,Q^2 \subset \bbR^r$ such that
    \begin{itemize}
        \item $\hat{z} \not\in Q^1$ and $\hat{z} \not\in Q^2$,
        \item $Q \supseteq Q^1 \cup Q^2$,
        \item $Q \cap H = (Q^1 \cap H) \cup (Q^2 \cap H)$, and
        \item $Q^1 \cap Q^2 = \emptyset$.
    \end{itemize}
\end{definition}

We note that branching is described solely in terms of the  control variables $z$, in contrast to a constraint branching approach \cite{Beale:1970,Farias-Jr.:2008,Farias-Jr.:2013,Keha:2006,Tomlin:1981}, which would work directly on the original variables $x$. Additionally, our branching schemes map back to our original setting in a straightforward way: if $R$ is the LP relaxation for our formulation in $(x,z)$-space, take $Q = \Proj_z(R)$. Then we can construct $Q^1$ and $Q^2$ by adding linear inequalities to $Q$. These inequalities will map to a linear inequality for $R$ with support only on the $z$ variables, giving two polyhedra $R^1$ and $R^2$ in the original $(x,z)$-space.

We will call $R^1$ and $R^2$ the LP relaxations for the \emph{subproblems}, and $Q^1$ and $Q^2$ the \emph{code relaxations} for the subproblems. In the case that an encoding $H$ has an associated branching scheme, we will say that the corresponding formulation is a \emph{mixed-integer branching formulation} to emphasize that this is a strict generalization of traditional MIP formulations, and that mixed-integer branching formulations retain many of the computational properties of traditional MIP formulations relevant for branch-and-bound and branch-and-cut methods.


Although a branch-and-bound method using variable branching may produce exponentially many subproblems, it enjoys a finite termination guarantee. This is not necessarily the case for any branching scheme constructed according to Definition~\ref{def:branching-scheme}. However, it is not difficult to see that, as $H$ is finite, a sufficient condition for finite termination is that $\Conv(Q^1 \cap H) = Q^1$ and $\Conv(Q^2 \cap H) = Q^2$.

\subsection{The reflected binary Gray and zig-zag encodings} \label{sec:brg-zigzag}

The first two encodings will will present have previously been used in the literature to construct traditional MIP formulations for the SOS2 constraint. Both are defined recursively~\cite{Huchette:2017} by the rows of the matrices $K^1 = C^1 = (0,1)^T$ and
\[
    K^{t+1} = \begin{pmatrix} K^t & {\bf 0}^t \\ \rev(K^t) & {\bf 1}^t \end{pmatrix}, \quad C^{t+1} = \begin{pmatrix} C^t & {\bf 0}^t \\ C^t + {\bf 1}^t \otimes C^t_t & {\bf 1}^t \end{pmatrix} \quad\quad \forall t \in \{2,3,\ldots\},
\]
where ${\bf 0}^t \in \bbR^t$ (respectively ${\bf 1}^t \in \bbR^t$) is the vector will all entries equal to zero (respectively one), $A_i$ is the $i$-th row of the matrix $A$, $\rev(A)$ reverses the rows of the matrix $A$, and $u \otimes v = u v^T \subseteq \bbR^{m \times n}$ for any $u \in \bbR^m$ and $v \in \bbR^n$.

\begin{figure}[htpb]
    \centering
    \begin{tikzpicture}
        [thick, ->] \draw (0,0,0) -- (1,0,0);
        [thick, ->] \draw (1,0,0) -- (1,1,0);
        [thick, ->] \draw (1,1,0) -- (0,1,0);
        [thick, ->] \draw (0,1,0) -- (0,1,1);
        [thick, ->] \draw (0,1,1) -- (1,1,1);
        [thick, ->] \draw (1,1,1) -- (1,0,1);
        [thick, ->] \draw (1,0,1) -- (0,0,1);
        \begin{scope}[canvas is xy plane at z=0]
            \draw[help lines] (0,0) grid (2,2);
            [ultra thick] \node [right] at (2,0) {{\smaller $z_1$}};
            [ultra thick] \node [left] at (0,0) {{\smaller $h^1$}};
        \end{scope}
        \begin{scope}[canvas is xy plane at z=1]
            \node [left] at (0,0) {\smaller $h^8$};
        \end{scope}
        \begin{scope}[canvas is xz plane at y=0]
            \draw[help lines] (0,0) grid (2,2);
            [ultra thick] \node [below] at (0,2) {{\smaller $z_3$}};
        \end{scope}
        \begin{scope}[canvas is yz plane at x=0]
            \draw[help lines] (0,0) grid (2,2);
            [ultra thick] \node [above] at (2,0) {{\smaller $z_2$}};
        \end{scope}
        
        \begin{scope}[canvas is xy plane at z=0]
            \draw [fill] (0,0) circle [radius=.05];
        \end{scope}
    \end{tikzpicture}
    \begin{tikzpicture}
        [thick, ->] \draw (0,0,0) -- (1,0,0);
        [thick, ->] \draw (1,0,0) -- (1,1,0);
        [thick, ->] \draw (1,1,0) -- (2,1,0);
        [thick, ->] \draw (2,1,0) -- (2,1,1);
        [thick, ->] \draw (2,1,1) -- (3,1,1);
        [thick, ->] \draw (3,1,1) -- (3,2,1);
        [thick, ->] \draw (3,2,1) -- (4,2,1);
        \begin{scope}[canvas is xy plane at z=0]
            \draw[help lines] (0,0) grid (5,2);
            [ultra thick] \node [right] at (5,0) {{\smaller $z_1$}};
            [ultra thick] \node [left] at (0,0) {{\smaller $h^1$}};
        \end{scope}
        \begin{scope}[canvas is xy plane at z=1]
            \node [right] at (4,2) {\smaller $h^8$};
        \end{scope}
        \begin{scope}[canvas is xz plane at y=0]
            \draw[help lines] (0,0) grid (5,2);
            [ultra thick] \node [below] at (0,2) {{\smaller $z_3$}};
        \end{scope}
        \begin{scope}[canvas is yz plane at x=0]
            \draw[help lines] (0,0) grid (2,2);
            [ultra thick] \node [above] at (2,0) {{\smaller $z_2$}};
        \end{scope}
        \begin{scope}[canvas is xy plane at z=0]
            \draw [fill] (0,0) circle [radius=.05];
        \end{scope}
    \end{tikzpicture}
    \caption{Depiction of the binary reflected Gray encoding $\HRB_8$ (\textbf{Left}) and the zig-zag encoding $\HZZ_8$ (\textbf{Right}).}
    \label{fig:hole-free}
\end{figure}
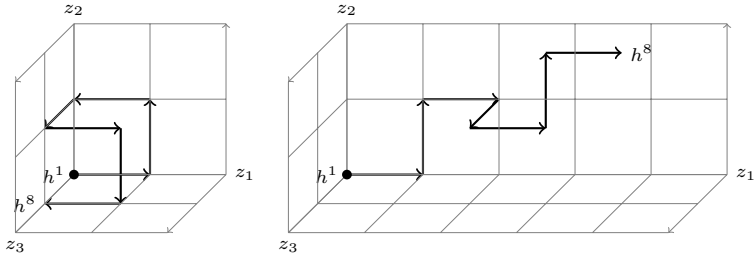

For some fixed number $d \in \bbN$, take $r = \lceil \log_2(d) \rceil$. We will take the \emph{reflected binary Gray} encoding~\cite{Huchette:2017,Savage:1997,Vielma:2010,Vielma:2009a} $\HRB_d \subseteq \{0,1\}^r$ as the first $d$ rows of $K^r$, and the \emph{zig-zag} encoding~\cite{Huchette:2017} $\HZZ_d$ as the first $d$ rows of $C^r$. In Figure~\ref{fig:hole-free}, we see the two encodings for $d=8$ (cf. Figure~1 of \cite{Huchette:2017}).

It can be shown that both the reflected binary Gray and zig-zag encodings are hole-free and in convex position~\cite{Huchette:2017}, and so therefore lead to traditional MIP formulations. Indeed, Huchette and Vielma~\cite{Huchette:2017} use both encodings to construct traditional MIP formulations for the SOS2 constraint that use $\lceil \log_2(d) \rceil$ control variables and $2\lceil \log_2(d) \rceil$ general inequality constraints. In Section~\ref{sec:MIP-annulus}, we will apply Theorem~\ref{thm:general-cdc-characterization} using the two encodings to construct logarithmic traditional MIP formulations for a relaxation of the annulus.

Additionally, the two hole-free encodings give us an opportunity to show how traditional variable branching fits into our branching scheme framework. Take either $H=\HRB_d$ or $H=\HZZ_d$, and consider some point $\hat{z} \in \Conv(H)$. If $\hat{z} \in \bbZ^r$, this verifies that $\hat{z} \in H$. Otherwise, we can select a component $k \in \llbracket r \rrbracket$ which is fractional, i.e. $\hat{z}_{k} \not\in \bbZ$. Then the two child code relaxations are created by rounding this component: $Q^1 = \{z \in Q : z_k \leq \lfloor \hat{z}_k \rfloor \}$ and $Q^2 = \{z \in Q : z_k \geq \lceil \hat{z}_k \rceil \}$.

\subsection{Moment curve encoding} \label{sec:moment-curve}

\begin{figure}[htpb]
    \centering
    \begin{tikzpicture}[xscale=0.8,yscale=0.2]
        \draw [dashed, fill=gray!20] (1,1) -- (2,4) -- (3,9) -- (4,16) -- (5,25) -- (6,36) -- (7,49) -- (1,1);

        \draw [fill=gray!80] (1,1) -- (2,4) -- (3,9) -- (1,1);
        \draw [fill=gray!80] (4,16) -- (5,25) -- (6,36) -- (7,49) -- (4,16);

        \draw [fill,yscale=4] (3.5,17.5/4) circle [radius=.05];

        \draw[help lines] (1,1) grid (7,49);
    \end{tikzpicture} \hspace{2em}
    \begin{tikzpicture}[xscale=0.8,yscale=0.2]
        \draw [dashed, fill=gray!20] (1,1) -- (2,4) -- (3,9) -- (4,16) -- (5,25) -- (6,36) -- (7,49) -- (1,1);
        \draw [fill=gray!80] (1,1) -- (2,4) -- (3,9) -- (4,16) -- (1,1);
        \draw [fill=gray!80]  -- (6,36) -- (7,49) -- (5,25);

        \draw [fill,yscale=4] (4,25/4) circle [radius=.05];

        \draw[help lines] (1,1) grid (7,49);
    \end{tikzpicture}

    \caption{Illustration of the branching scheme for the moment curve encoding $\HMC_7$. The original code relaxation in the $z$-space is shown in the dashed region, and those for the two subproblems are shown in the darker shaded regions. The optimal solution for the original LP relaxation is depicted with a solid dot. We show the branching with a solution that is fractional \textbf{(Left)}, and one where there is no valid variable branching disjunction to separate the point \textbf{(Right)}.}
    \label{fig:moment-curve-branching}
\end{figure}
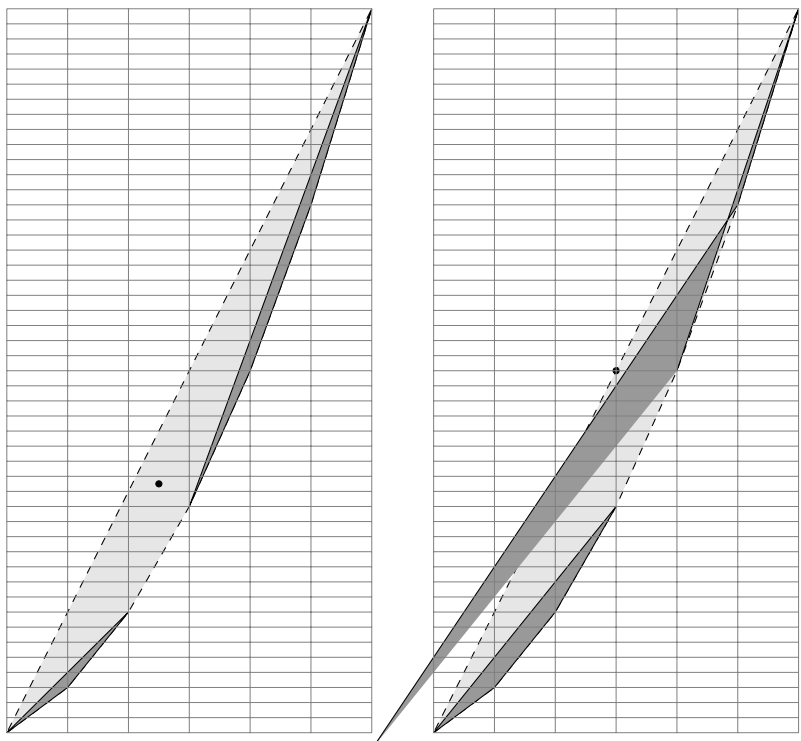

The $\eta$-dimensional moment curve is given by the function $m_\eta(t) = (t,t^2,\ldots,t^\eta)$. Given $d (\geq \eta)$ ordered points $t_1 < t_2 < \cdots < t_d$ on the real line, the corresponding \emph{cyclic polytope} is $\Conv(\{m_\eta(t_i)\}_{i=1}^d)$, a well-studied object~\cite{Bogomolov:2015,Ziegler:2007}. For our purposes, we are interested in constructing encodings that lie along the two-dimensional moment curve: $\HMC_d \defeq (m_2(i))_{i=1}^d$. If $d > 2$, then this choice of encoding is not hole-free; for example, $\frac{1}{2}(m_2(1)+m_2(3)) = (2,5) \not\in \HMC_d$. However, the encoding is in convex position, and it is straightforward to check if $\hat{z} \in \HMC_d$. We also see that a description for $\Psi_d(l,u) \defeq \Conv(\{z \in \HMC_d :  l \leq z_1 \leq u\})$ is
\begin{subequations} \label{eqn:moment-curve-convex-hull}
\begin{align}
    z_2 - i^2 &\geq (2i+1)(z_1-i) \quad \forall i \in \llbracket l,u-1 \rrbracket \\
    (u-l)(z_2 - l^2) &\leq (u^2-l^2)(z_1-l).
\end{align}
\end{subequations}
Our branching scheme for the encoding $\HMC_d$ starts with a relaxation of the form $Q = \Psi_d(\ell,u)$ for some $\ell, u \in \bbZ$. Provided that $\hat{z} \not\in \HMC_d$, we create two child code relaxations of the form $Q^1 = \Psi_d(\ell,\lfloor \hat{z}_1 \rfloor)$ and $Q^2 = \Psi_d(\lfloor \hat{z}_1 \rfloor + 1,u)$. See Figure~\ref{fig:moment-curve-branching} for an illustration of the branching.

We emphasize that while this branching scheme uses two-term disjunction branching, in nearly every case a (potentially wide) variable branching disjunction is also valid. For example, the variable branching disjunction $z_1 \leq 2 \vee z_1 \geq 3$ is valid for the point in the left side of Figure~\ref{fig:moment-curve-branching}. This will be the case in all but pathological cases: for example, that depicted in the right side of Figure~\ref{fig:moment-curve-branching}. This means that the branching portion of the algorithm can proceed using the branching scheme described above, while the cut generation procedure can also use the valid variable branching split disjunctions as well.

\subsection{A more exotic encoding} \label{sec:exotic-encoding}
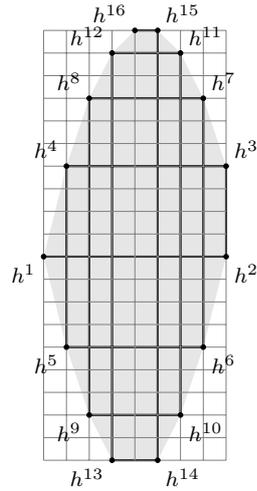
\begin{figure}[htpb]
    \centering
    \begin{tikzpicture}[scale=.3]
        \fill [fill=gray!20] (0,0) -- (1,4) -- (2,7) -- (3,9) -- (4,10) -- (5,10) -- (6,9) -- (7,7) -- (8,4) -- (8,0) -- (7,-4) -- (6,-7) -- (5,-9) -- (3,-9) -- (2,-7) -- (1,-4) -- (0,0);

        \draw [thick] (0,0) -- (8,0) -- (8,4) -- (1,4) -- (1,-4) -- (7,-4) -- (7,7) -- (2,7) -- (2,-7) -- (6,-7) -- (6,9) -- (3,9) -- (3,-9) -- (5,-9) -- (5,10) -- (4,10);

        \draw[help lines] (0,-9) grid (8,10);
         [ultra thick] \node [below left] at (0,0) {$h^1$};
         \draw [fill] (0,0) circle [radius=.1];

         [ultra thick] \node [below right] at (8,0) {$h^2$};
         \draw [fill] (8,0) circle [radius=.1];

         [ultra thick] \node [above right] at (8,4) {$h^3$};
         \draw [fill] (8,4) circle [radius=.1];

         [ultra thick] \node [above left] at (1,4) {$h^4$};
         \draw [fill] (1,4) circle [radius=.1];

         [ultra thick] \node [below left] at (1,-4) {$h^5$};
         \draw [fill] (1,-4) circle [radius=.1];

         [ultra thick] \node [below right] at (7,-4) {$h^6$};
         \draw [fill] (7,-4) circle [radius=.1];

         [ultra thick] \node [above right] at (7,7) {$h^7$};
         \draw [fill] (7,7) circle [radius=.1];

         [ultra thick] \node [above left] at (2,7) {$h^8$};
         \draw [fill] (2,7) circle [radius=.1];

         [ultra thick] \node [below left] at (2,-7) {$h^9$};
         \draw [fill] (2,-7) circle [radius=.1];

         [ultra thick] \node [below right] at (6,-7) {$h^{10}$};
         \draw [fill] (6,-7) circle [radius=.1];

         [ultra thick] \node [above right] at (6,9) {$h^{11}$};
         \draw [fill] (6,9) circle [radius=.1];

         [ultra thick] \node [above left] at (3,9) {$h^{12}$};
         \draw [fill] (3,9) circle [radius=.1];

         [ultra thick] \node [below left] at (3,-9) {$h^{13}$};
         \draw [fill] (3,-9) circle [radius=.1];

         [ultra thick] \node [below right] at (5,-9) {$h^{14}$};
         \draw [fill] (5,-9) circle [radius=.1];

         [ultra thick] \node [above right] at (5,10) {$h^{15}$};
         \draw [fill] (5,10) circle [radius=.1];

         [ultra thick] \node [above left] at (4,10) {$h^{16}$};
         \draw [fill] (4,10) circle [radius=.1];

    \end{tikzpicture}
    \caption{The exotic two-dimensional encoding $\HEX_{16}$.}
    \label{fig:sos2-constant}
\end{figure}

Consider the two-dimensional encoding in Figure~\ref{fig:sos2-constant}. Take some positive integer $r$, along with $d=4r$, and consider the encoding $\HEX_d = (h^i)_{i=1}^{d}$ where
\begin{subequations} \label{eqn:exotic-encoding}
\begin{align}
    h^{4k-3} &= \left( k-r-1,  \frac{1}{2}(k-1)(k-2r-2) \right) \\
    h^{4k-2} &= \left( r-k+1, \frac{1}{2}(k-1)(k-2r-2) \right) \\
    h^{4k-1} &= \left( r-k+1,  -\frac{1}{2}k(k-2r-1)  \right) \\
    h^{4k}   &= \left( k-r,   -\frac{1}{2}k(k-2r-1)  \right)
\end{align}
\end{subequations}
for each $k \in \llbracket r \rrbracket$. These points are in convex position.

\begin{proposition}
    For any $r \in \bbN$, the points $\HEX_{4r}$ are in convex position.
\end{proposition}
\proof{}
    The result for $r=1$ follows from inspection, so presume that $r > 1$. For each point $h^i$, we propose an inequality $c^i \cdot z \leq b^i$ that strictly separates $h^i$ from the remaining codes in $\HEX_d$. For each $k \in \llbracket r \rrbracket$, the coefficients are
    \begin{alignat*}{2}
        c^{4k-3} &= \left( -(r-k+2)-(r-k+1), -2 \right) \\
        c^{4k-2} &= \left(  (r-k+2)+(r-k+1), -2 \right) \\
        c^{4k-1} &= \left(  (r-k+1)+(r-k), 2 \right) \\
        c^{4k}   &= \left( -(r-k+1)-(r-k), 2 \right),
    \end{alignat*}
    where $b^i = c^{i} \cdot h^{i+4}$ for $i \in \llbracket 4 \rrbracket$ and $b^{i} = c^{i} \cdot h^{i-4}$ for $i \in \llbracket 5,4r \rrbracket$.
\qed \endproof

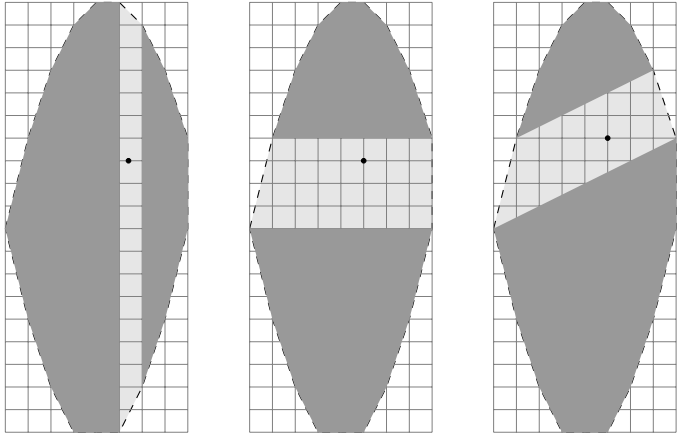
\begin{figure}
    \centering
    \begin{tikzpicture}[scale=.3]
        \draw [dashed, fill=gray!20] (0,0) -- (1,4) -- (2,7) -- (3,9) -- (4,10) -- (5,10) -- (6,9) -- (7,7) -- (8,4) -- (8,0) -- (7,-4) -- (6,-7) -- (5,-9) -- (3,-9) -- (2,-7) -- (1,-4) -- (0,0);
        
        \draw[help lines] (0,-9) grid (8,10);

        \fill [fill=gray!80] (0,0) -- (1,4) -- (2,7) -- (3,9) -- (4,10) -- (5,10) -- (5,-9) -- (3,-9) -- (2,-7) -- (1,-4) -- (0,0);

        \fill [fill=gray!80] (6,9) -- (7,7) -- (8,4) -- (8,0) -- (7,-4) -- (6,-7);

        \draw [fill] (5.4,3) circle [radius=0.1];

    \end{tikzpicture} \hspace{2em}
    \begin{tikzpicture}[scale=.3]
        \draw [dashed, fill=gray!20] (0,0) -- (1,4) -- (2,7) -- (3,9) -- (4,10) -- (5,10) -- (6,9) -- (7,7) -- (8,4) -- (8,0) -- (7,-4) -- (6,-7) -- (5,-9) -- (3,-9) -- (2,-7) -- (1,-4) -- (0,0);

        \draw[help lines] (0,-9) grid (8,10);

        \fill [fill=gray!80] (1,4) -- (2,7) -- (3,9) -- (4,10) -- (5,10) -- (6,9) -- (7,7) -- (8,4);

        \fill [fill=gray!80] (0,0) -- (8,0) -- (7,-4) -- (6,-7) -- (5,-9) -- (3,-9) -- (2,-7) -- (1,-4) -- (0,0);

        \draw [fill] (5,3) circle [radius=0.1];

    \end{tikzpicture} \hspace{2em}
    \begin{tikzpicture}[scale=.3]
        \draw [dashed, fill=gray!20] (0,0) -- (1,4) -- (2,7) -- (3,9) -- (4,10) -- (5,10) -- (6,9) -- (7,7) -- (8,4) -- (8,0) -- (7,-4) -- (6,-7) -- (5,-9) -- (3,-9) -- (2,-7) -- (1,-4) -- (0,0);

        \draw[help lines] (0,-9) grid (8,10);

        \fill [fill=gray!80] (1,4) -- (2,7) -- (3,9) -- (4,10) -- (5,10) -- (6,9) -- (7,7);

        \fill [fill=gray!80] (0,0) -- (8,4) -- (8,0) -- (7,-4) -- (6,-7) -- (5,-9) -- (3,-9) -- (2,-7) -- (1,-4) -- (0,0);

        \draw [fill] (5,4) circle [radius=0.1];
    \end{tikzpicture}
    \caption{Branching scheme for the exotic encoding $\HMC_{16}$ ($r=4$) when the LP optimal solution for the control variables $\hat{z}$ has: \textbf{(Left)} $\hat{z}_1$ fractional, \textbf{(Center)} $\hat{z}_1 \in \bbZ$ but $\hat{z}_2 \notin Y = \{h_2 : h \in \HMC_{16}\}$, and \textbf{(Right)} $\hat{z}_1 \in \bbZ$, $\hat{z}_2 \in Y$, and $\hat{z} \notin \HMC_{16}$. The relaxations for the two subproblems in each are the two shaded regions in each picture.}
    \label{fig:sos2-constant-branching}
\end{figure}

The structure of this encoding also suggests a relatively simple branching scheme. Given a point $\hat{z} \not\in \HEX_d$, we consider three cases, depicted in Figure~\ref{fig:sos2-constant-branching}. In the first case, $\hat{z}_1 \not\in \bbZ$, and we perform standard variable branching: $Q^1 = \{z \in Q : z_1 \leq \lfloor \hat{z}_1 \rfloor\}$ and $Q^2 = \{z \in Q : z_1 \geq \lceil \hat{z}_1 \rceil\}$. If $\hat{z}_1 \in \bbZ$, then we consider two other cases. Take $Y = \{h_2 : h \in \HEX_d\}$ as the set of all values the encoding takes in the second component, $\underline{b} = \max\{t \in Y : t < \hat{z}_2\}$, and $\overline{b} = \min\{t \in Y : t > \hat{z}_2\}$. If $\hat{z}_2 \not \in Y$, then we apply a wide variable branching of the form $Q^1 = \{z \in Q : z_2 \leq \underline{b}\}$, and $Q^2 = \{z \in Q : z_2 \geq \overline{b}\}$.

The final case remains where $\hat{z}_1 \in \bbZ$, $\hat{z}_2 \in Y$, and yet $\hat{z} \not\in \HEX_d$. In this case, we will branch on a two-term non-parallel disjunction. Take $Z(b) = \{h \in \HEX_d : h_2 = b\}$. We take the nearest point to the northeast of $\hat{z}$ as $h^{NE} \in Z(\overline{b})$ such that $h^{NE}_1 \geq \max_{h \in Z(\overline{b})}h_1$. Next, take the nearest point to the southwest $h^{SW} \in Z(\underline{b})$ such that $h^{SW} \leq \max_{h \in Z(\underline{b})}h_1$. Take the points directly to the west and east of $\hat{z}$, $h^{W}, h^{E} \in Z(\hat{z}_2)$ (i.e. $h^W_1 < h^E_1$), and we can express the two-term non-parallel disjunction branching with two child code relaxations as $Q^1 = \{z \in Q : (h^{NE}_1-h^W_1)(z_2-h^W_2) \geq (h^{NE}_2-h^W_2)(z_1-h^W_1)\}$ and $Q^2 = \{z \in Q : (h^{SW}_1-h^E_1)(z_2-h^E_2) \geq (h^{SW}_2-h^E_2)(z_1-h^E_1)\}$.

\section{Very small mixed-integer branching formulations}

We are now in a position to derive very small mixed-integer branching formulations for combinatorial disjunctive constraints. Each formulation will have only two control variables, and will be constructed using the two-dimensional encodings presented in the previous section. Along the way, we  also present two new logarithmic-sized traditional MIP formulations for a relaxation of the annulus that follow as a natural consequence of Theorem~\ref{thm:general-cdc-characterization}. Combined, these results illustrate that Theorem~\ref{thm:general-cdc-characterization} can be practically used to construct both traditional MIP and mixed-integer branching formulations for disjunctive constraints. 

\subsection{Very small formulations for general combinatorial disjunctive constraints}




First, we state a general result: given any combinatorial disjunctive constraint and any two-dimensional encoding in convex position, we can provide an explicit description for a very small ideal formulation.

\begin{proposition} \label{prop:general-cdc-with-general-2D-embedding}
    Take $\calT = (T^i \subseteq \llbracket n \rrbracket)_{i=1}^d$ and let $H = (h^i)_{i=1}^d \subset \bbR^2$ be a two-dimensional encoding in convex position. Take $b^{i,j} = (c^{i,j}_2,-c^{i,j}_1)$ for each $\{i,j\} \in [d]^2$. Then $(\lambda,z) \in Q(\calP(\calT),H)$ if and only if
    \begin{gather*}
        \sum_{v=1}^n \min_{s : v \in T^s} \{b^{i,j} \cdot h^s\} \lambda_v \leq b^{i,j} \cdot z \leq \sum_{v=1}^n \max_{s : v \in T^s} \{b^{i,j} \cdot h^s\} \lambda_v \quad \forall \{i,j\} \in [d]^2 \\
        (\lambda,z) \in \Delta^{n} \times \bbR^2.
    \end{gather*}
\end{proposition}
\proof{}
    The result follows from Theorem~\ref{thm:general-cdc-characterization}. If $D$ is not connected, we may introduce an artificial $\lambda_{n+1}$ variable to the constraint, and append it $T \leftarrow T \cup \{n+1\}$ to each set $T \in \calT$. The corresponding edge set $D' = [d]^2$ is now connected, and we can simply impose that $\lambda_{n+1} \leq 0$ to recover our original constraint.

    First, we observe that $b^{i,j} \cdot c^{i,j} = 0$, and so as $\calL$ is two-dimensional, $M(b^{i,j};\calL)$ is the hyperplane spanned by $c^{i,j}$. Furthermore, we have that $D = \{\{i,j\} \in [d]^2 : T^i \cap T^j \neq \emptyset\} \subseteq D' = [d]^2$, and so this representation will recover all the inequalities in \eqref{eqn:general-V-formulation-1}. It just remains to show that any inequality given by $\{i,j\} \in D' \backslash D$ is valid for $Q(\calP(\calT),H)$. To see this, consider any $(\lambda,z) = ({\bf e}^w,h^u) \in \Em(\calP(\calT),H)$; that is, $w \in T^u$. Then $\sum_{v=1}^n \min_{s : v \in T^s}\{b^{i,j} \cdot h^s\} \lambda_v = \min_{s : w \in T^s}\{b^{i,j} \cdot h^s\} \leq b^{i,j} \cdot h^u$, as $b^{i,j} \cdot h^u$ is one of the terms appearing in the minimization. A similar argument holds for the other side of the constraint.
\qed \endproof

This result implies a quadratic $\scrO(d^2)$ upper bound on the number of general inequality constraints needed to construct an ideal mixed-integer branching formulations for \emph{any} combinatorial disjunctive constraint. This is in sharp contrast to the traditional MIP setting, where binary encodings can--and typically do--lead to an exponential number of facets~\cite{Vielma:2015a}.

Furthermore, this can be strengthened to an $\scrO(d)$ upper bound on the number of general inequality constraints when we use the moment curve encoding.

\begin{proposition} \label{prop:general-cdc-moment-curve}
    Take $\calT = (T^i \subseteq \llbracket n \rrbracket)_{i=1}^d$. Then $(\lambda,z) \in Q(\calP(\calT),\HMC_d)$ if and only if
    \begin{gather*}
        \sum_{v=1}^n \min_{s : v \in T^s}\{s(t-s)\} \lambda_v \leq tz_1 - z_2 \leq \sum_{v=1}^n \max_{s : v \in T^s} \{s(t-s)\} \lambda_v \quad \forall t \in \llbracket 3, 2d-1 \rrbracket \\
        (\lambda,z) \in \Delta^{n} \times \bbR^2.
    \end{gather*}
\end{proposition}
\proof{}
    Take any $\{i,j\} \in [d]^2$. Observe that $c^{i,j} \equiv h^j-h^i = (j-i,j^2-i^2) = (j-i)\cdot(1,i+j)$, and that $3 \leq i+j \leq 2d-1$. Therefore, for each $\{i,j\} \in [d]^2$, there is some $t \in \llbracket 3,2d-1 \rrbracket$ and some $\alpha > 0$ such that $c^{i,j} = \alpha \cdot (1,t)$. Therefore, our representation here is equivalent to that in Proposition~\ref{prop:general-cdc-with-general-2D-embedding}, up to constant nonzero scalings of some of the inequalities.
\qed \endproof

\begin{figure}[htpb]
    \centering
    \begin{tikzpicture}
        \draw [fill=gray!20] (0,0) -- (1,0) -- (0,1) -- (0,0);
        \draw [fill=gray!20] (1,0) -- (1,1) -- (0,1) -- (1,0);
        \draw [fill=gray!20] (1,0) -- (2,0) -- (1,1) -- (1,0);
        \draw [fill=gray!20] (2,0) -- (2,1) -- (1,1) -- (2,0);
        \draw [fill=gray!20] (0,1) -- (1,1) -- (0,2) -- (0,1);
        \draw [fill=gray!20] (1,1) -- (1,2) -- (0,2) -- (1,1);
        \draw [fill=gray!20] (1,1) -- (2,1) -- (1,2) -- (1,1);
        \draw [fill=gray!20] (2,1) -- (2,2) -- (1,2) -- (2,1);

        \node [below left] at (0,0) {$1$};
        \node [below right] at (1,0) {$2$};
        \node [below right] at (2,0) {$3$};
        \node [below left] at (0,1) {$4$};
        \node [above right] at (1,1) {$5$};
        \node [below right] at (2,1) {$6$};
        \node [above left] at (0,2) {$7$};
        \node [above right] at (1,2) {$8$};
        \node [above right] at (2,2) {$9$};

        \draw [fill] (0,0) circle [radius=.05];
        \draw [fill] (1,0) circle [radius=.05];
        \draw [fill] (2,0) circle [radius=.05];
        \draw [fill] (0,1) circle [radius=.05];
        \draw [fill] (1,1) circle [radius=.05];
        \draw [fill] (2,1) circle [radius=.05];
        \draw [fill] (0,2) circle [radius=.05];
        \draw [fill] (1,2) circle [radius=.05];
        \draw [fill] (2,2) circle [radius=.05];
    \end{tikzpicture}
    \caption{A grid triangulation on the plane with $8$ alternatives (triangles). The nodes, or vertices for the triangles, are numbered.}
    \label{fig:triangulation}
\end{figure}
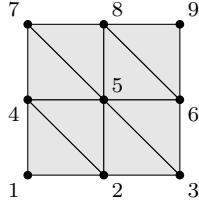

As a concrete example, consider the grid triangulation in Figure~\ref{fig:triangulation}. The sets $\calT = (T^i)_{i=1}^8$ correspond to each of the triangles, where
\begin{alignat*}{4}
    T^1 = \{1,2,4\}, \quad T^2 = \{5,6,8\}, \quad T^3 = \{3,5,6\}, \quad T^4 = \{4,5,7\}, \\ T^5 = \{5,7,8\}, \quad T^6 = \{2,3,5\}, \quad T^7 = \{2,4,5\}, \quad T^8 = \{6,8,9\}.
\end{alignat*}
Then a description for $Q(\calP(\calT),\HMC_8)$ using the moment curve encoding is:
\begin{align*}
    4\lambda_1 + 4\lambda_2 + 6\lambda_3 + 4\lambda_4 + 6\lambda_5 + 6\lambda_6 + 4\lambda_7 + 6\lambda_8 - 24\lambda_9 &\geq 5z_1 - z_2 \\
    6\lambda_1 + 6\lambda_2 + 12\lambda_3 + 12\lambda_4 + 12\lambda_5 + 12\lambda_6 + 12\lambda_7 + 10\lambda_8 - 8\lambda_9 &\geq 7z_1 - z_2 \\
    7\lambda_1 + 7\lambda_2 + 12\lambda_3 + 7\lambda_4 + 7\lambda_5 + 0\lambda_6 + 15\lambda_7 + 0\lambda_8 + 0\lambda_9 &\leq 8z_1 - z_2 \\
    8\lambda_1 + 8\lambda_2 + 18\lambda_3 + 8\lambda_4 + 14\lambda_5 + 8\lambda_6 + 20\lambda_7 + 8\lambda_8 + 8\lambda_9 &\leq 9z_1 - z_2 \\
    8\lambda_1 + 18\lambda_2 + 18\lambda_3 + 20\lambda_4 + 20\lambda_5 + 18\lambda_6 + 20\lambda_7 + 20\lambda_8 + 8\lambda_9 &\geq 9z_1 - z_2 \\
    9\lambda_1 + 9\lambda_2 + 21\lambda_3 + 9\lambda_4 + 16\lambda_5 + 16\lambda_6 + 24\lambda_7 + 16\lambda_8 + 16\lambda_9 &\leq 10z_1 - z_2 \\
    10\lambda_1 + 30\lambda_2 + 30\lambda_3 + 28\lambda_4 + 30\lambda_5 + 24\lambda_6 + 30\lambda_7 + 30\lambda_8 + 24\lambda_9 &\geq 11z_1 - z_2 \\
    12\lambda_1 + 42\lambda_2 + 42\lambda_3 + 42\lambda_4 + 42\lambda_5 + 40\lambda_6 + 40\lambda_7 + 40\lambda_8 + 40\lambda_9 &\geq 13z_1 - z_2 \\
    (\lambda,z) \in \Delta^9 \times \bbR^2.
\end{align*}
The construction of Proposition~\ref{prop:general-cdc-moment-curve} gives these 8 facet-inducing general inequality constraints, along with 16 others that are valid but not facet-inducing for $Q(\calP(\calT),\HMC_8)$, and therefore are not necessary. In contrast, any ideal binary MIP formulation (i.e. the encoding $H$ is some ordering of $\{0,1\}^3$) requires three control variables and at least 9 general inequality constraints~\cite{Huchette:2017}.

\subsection{A very small formulation for the SOS2 constraint}

We can sharpen our general results from the previous subsection if we take advantage of structure and choose an encoding tailored for a particular constraint. For example, the exotic encoding $\HEX_d$ was specifically designed for the SOS2 constraint, which we recall is given by the sets $\calT^{\SOS}_d = (\{i,i+1\})_{i=1}^d$.

\begin{proposition}\label{prop:sos2-constant}
    Take $d=4r$ for some $r \in \bbN$, and label $\HEX_d = (h^i)_{i=1}^d \subset \bbR^2$. Then $(\lambda,y) \in Q(\calP(\calT^{\SOS}_d),\HEX_d)$ if and only if
    \begin{subequations} \label{eqn:sos2-constant}
    \begin{align}
        h^1_k \lambda_1 + \sum_{i=2}^{d} \min\{h^{i-1}_k,h^{i}_k\}\lambda_i + h^{d}_k\lambda_{d+1} &\leq z_k \quad\quad \forall k \in \llbracket 2 \rrbracket \\
        h^1_k \lambda_1 + \sum_{i=2}^{d} \max\{h^{i-1}_k,h^{i}_k\}\lambda_i + h^{d}_k\lambda_{d+1} &\geq z_k \quad\quad \forall k \in \llbracket 2 \rrbracket \\
        (\lambda,z) \in \Delta^{d+1} \times \bbR^{2}.
    \end{align}
    \end{subequations}
\end{proposition}
\proof{}
    Apply Theorem~\ref{thm:general-cdc-characterization}, after observing that $C = \{c^{i,i+1} \equiv h^{i+1}-h^{i}\}_{i=1}^{d-1} \subseteq \{\pm {\bf e}^1, \pm {\bf e}^2\}$, and so taking $b^1 = {\bf e}^1$ and $b^2={\bf e}^2$ suffices.
\qed \endproof

As a concrete example, the formulation \eqref{eqn:sos2-constant} for the SOS2 constraint on $n=17$ components is
\begin{subequations}
\begin{gather}
    \notag -4\lambda_1 -4\lambda_{2} + 4\lambda_{3} -3\lambda_{4} -3\lambda_{5} -3\lambda_{6} + 3\lambda_{7} -2\lambda_{8} -2\lambda_{9} + \\ -2\lambda_{10} + 2\lambda_{11} -1\lambda_{12} -1\lambda_{13} -1\lambda_{14} + 1\lambda_{15} + 0\lambda_{16} + 0\lambda_{17} \leq z_1 \\
    \notag -4\lambda_1 + 4\lambda_{2} + 4\lambda_{3} + 4\lambda_{4} -3\lambda_{5} + 3\lambda_{6} + 3\lambda_{7} + 3\lambda_{8} -2\lambda_{9} + \\ 2\lambda_{10} + 2\lambda_{11} + 2\lambda_{12} -1\lambda_{13} + 1\lambda_{14} + 1\lambda_{15} + 1\lambda_{16} + 0\lambda_{17} \geq z_1 \\
    \notag 0\lambda_1 + 0\lambda_{2} + 0\lambda_{3} + 4\lambda_{4} -4\lambda_{5} -4\lambda_{6} -4\lambda_{7} + 7\lambda_{8} -7\lambda_{9} + \\ -7\lambda_{10} -7\lambda_{11} + 9\lambda_{12} -9\lambda_{13} -9\lambda_{14} -9\lambda_{15} + 10\lambda_{16} + 10\lambda_{17} \leq z_2 \\
    \notag 0\lambda_1 + 0\lambda_{2} + 4\lambda_{3} + 4\lambda_{4} + 4\lambda_{5} -4\lambda_{6} + 7\lambda_{7} + 7\lambda_{8} + 7\lambda_{9} + \\ -7\lambda_{10} + 9\lambda_{11} + 9\lambda_{12} + 9\lambda_{13} -9\lambda_{14} + 10\lambda_{15} + 10\lambda_{16} + 10\lambda_{17} \leq z_2 \\
    (\lambda,z) \in \Delta^{17} \times \bbR^2.
\end{gather}
\end{subequations}
This ideal mixed-integer branching formulation uses only two control variables, along with only four general integer inequality constraints. We contrast this with the logarithmic formulations of Huchette and Vielma~\cite{Huchette:2017,Vielma:2010,Vielma:2009a}, which are also ideal but require 4 control variables and 8 general inequality constraints. Moreover, the numbers of control variables or general inequality constraints in formulation \eqref{eqn:sos2-constant} do not grow with $d$, and so this difference will be even more pronounced with larger instances of the SOS2 constraint.

\subsection{Relaxations of the annulus}

\begin{figure}[htpb]
    \centering
    \begin{tikzpicture}[scale=.7]

        \draw [<->] (-4,0) -- (4,0);
        \draw [<->] (0,-4) -- (0,4);
        \node [left] at (0,4) {$x_2$};
        \node [above] at (4,0) {$x_1$};
        \draw[fill=gray!50,even odd rule]  circle (3) circle (2);

        \draw [<->,dashed] (0,0) -- (1.7320508075688774,1);
        \draw [<->,dashed] (0,0) -- (1.5,2.5980762113533156);

        \node [below right] at (1.7320508075688774/2,1/2) {$s$};
        \node [above left] at (1.5/2,2.5980762113533156/2) {$S$};

    \begin{scope}[shift={(9,0)}]
        \draw [fill=gray!30] (1.4142135623730951,1.414213562373095) -- (2.0,0.0) -- (3.2471766008771823,0.0) -- (2.296100594190539,2.2961005941905386);
        \draw [fill=gray!30] (1.2246467991473532e-16,2.0) -- (1.4142135623730951,1.414213562373095) -- (2.296100594190539,2.2961005941905386) -- (1.988322215265212e-16,3.2471766008771823);
        \draw [fill=gray!30] (-1.414213562373095,1.4142135623730951) -- (1.2246467991473532e-16,2.0) -- (1.988322215265212e-16,3.2471766008771823) -- (-2.2961005941905386,2.296100594190539);
        \draw [fill=gray!30] (-2.0,2.4492935982947064e-16) -- (-1.414213562373095,1.4142135623730951) -- (-2.2961005941905386,2.296100594190539) -- (-3.2471766008771823,3.976644430530424e-16);
        \draw [fill=gray!30] (-1.4142135623730954,-1.414213562373095) -- (-2.0,2.4492935982947064e-16) -- (-3.2471766008771823,3.976644430530424e-16) -- (-2.2961005941905395,-2.2961005941905386);
        \draw [fill=gray!30] (-3.6739403974420594e-16,-2.0) -- (-1.4142135623730954,-1.414213562373095) -- (-2.2961005941905395,-2.2961005941905386) -- (-5.964966645795635e-16,-3.2471766008771823);
        \draw [fill=gray!30] (1.4142135623730947,-1.4142135623730954) -- (-3.6739403974420594e-16,-2.0) -- (-5.964966645795635e-16,-3.2471766008771823) -- (2.296100594190538,-2.2961005941905395);
        \draw [fill=gray!30] (2.0,-4.898587196589413e-16) -- (1.4142135623730947,-1.4142135623730954) -- (2.296100594190538,-2.2961005941905395) -- (3.2471766008771823,-7.953288861060848e-16);
        \draw[fill=gray!50,even odd rule,dashed]  circle (3) circle (2);

        \node at (2.309698831278217, 0.9567085809127245) {$P^1$};
        \node at (0.9567085809127246, 2.309698831278217) {$P^2$};
        \node at (-0.9567085809127243, 2.309698831278217) {$P^3$};
        \node at (-2.309698831278217, 0.9567085809127247) {$P^4$};
        \node at (-2.309698831278217, -0.9567085809127241) {$P^5$};
        \node at (-0.9567085809127258, -2.309698831278216) {$P^6$};
        \node at (0.956708580912725, -2.3096988312782165) {$P^7$};
        \node at (2.309698831278216, -0.956708580912726) {$P^8$};

        \draw (1.4142135623730951,1.414213562373095) -- (2.0,0.0) -- (3.2471766008771823,0.0) -- (2.296100594190539,2.2961005941905386);
        \draw (1.2246467991473532e-16,2.0) -- (1.4142135623730951,1.414213562373095) -- (2.296100594190539,2.2961005941905386) -- (1.988322215265212e-16,3.2471766008771823);
        \draw (-1.414213562373095,1.4142135623730951) -- (1.2246467991473532e-16,2.0) -- (1.988322215265212e-16,3.2471766008771823) -- (-2.2961005941905386,2.296100594190539);
        \draw (-2.0,2.4492935982947064e-16) -- (-1.414213562373095,1.4142135623730951) -- (-2.2961005941905386,2.296100594190539) -- (-3.2471766008771823,3.976644430530424e-16);
        \draw (-1.4142135623730954,-1.414213562373095) -- (-2.0,2.4492935982947064e-16) -- (-3.2471766008771823,3.976644430530424e-16) -- (-2.2961005941905395,-2.2961005941905386);
        \draw (-3.6739403974420594e-16,-2.0) -- (-1.4142135623730954,-1.414213562373095) -- (-2.2961005941905395,-2.2961005941905386) -- (-5.964966645795635e-16,-3.2471766008771823);
        \draw (1.4142135623730947,-1.4142135623730954) -- (-3.6739403974420594e-16,-2.0) -- (-5.964966645795635e-16,-3.2471766008771823) -- (2.296100594190538,-2.2961005941905395);
        \draw (2.0,-4.898587196589413e-16) -- (1.4142135623730947,-1.4142135623730954) -- (2.296100594190538,-2.2961005941905395) -- (3.2471766008771823,-7.953288861060848e-16);

        \node [left] at (2.0,0.0) {$v^{-1} \equiv v^{15}$};
        \node [right] at (3.2471766008771823,0.0) {$v^{0} \equiv v^{16}$};
        \node [below left] at (1.4142135623730951,1.414213562373095) {$v^{1}$};
        \node [above right] at (2.296100594190539,2.2961005941905386) {$v^{2}$};
        \node [below] at (1.2246467991473532e-16,2.0) {$v^{3}$};
        \node [above] at (1.988322215265212e-16,3.2471766008771823) {$v^{4}$};
        \node [below right] at (-1.414213562373095,1.4142135623730951) {$v^{5}$};
        \node [above left] at (-2.2961005941905386,2.296100594190539) {$v^{6}$};
        \node [right] at (-2,0) {$v^{7}$};
        \node [left] at (-3.2471766008771823,0.0) {$v^{8}$};
        \node [above right] at (-1.414213562373095,-1.4142135623730951) {$v^{9}$};
        \node [below left] at (-2.2961005941905386,-2.296100594190539) {$v^{10}$};
        \node [above] at (-1.2246467991473532e-16,-2.0) {$v^{11}$};
        \node [below] at (-1.988322215265212e-16,-3.2471766008771823) {$v^{12}$};
        \node [above left] at (1.4142135623730951,-1.414213562373095) {$v^{13}$};
        \node [below right] at (2.296100594190539,-2.2961005941905386) {$v^{14}$};

    \end{scope}
    \end{tikzpicture}
    \caption{(\textbf{Left}) The annulus $\calA$ and (\textbf{Right}) its corresponding quadrilateral relaxation $\hat{\calA}$ given by \eqref{eqn:annulus-pieces} with $d=8$.}
    \label{fig:annulus}
\end{figure}
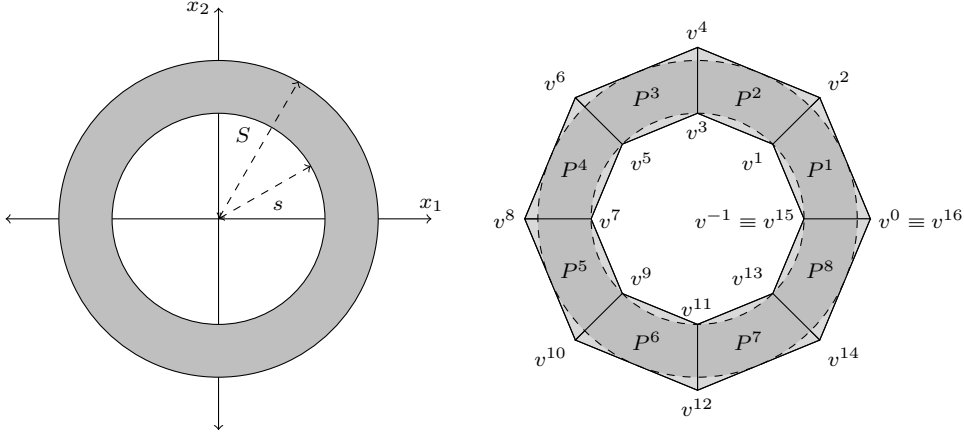

The annulus is a set in the plane $\calA = \{x \in \bbR^2 : s \leq ||x||_2 \leq S\}$ for constants $s, S \in \bbR_{+}$; see the left side of Figure~\ref{fig:annulus} for an illustration. A constraint of the form $x \in \calA$ might arise when modeling a complex number $z = x_1 + x_2 \mathbf{i}$, as $x \in \calA$ bounds the magnitude of $z$ as $s \leq |z| \leq S$. Such constraints arise in power systems optimization: for example, in the ``rectangular formulation''~\cite{Kocuk:2015} and the second-order cone reformulation~\cite{Jabr:2006,Liu:2017} of the optimal power flow problem, and the reactive power dispatch problem~\cite{Foster:2013}. Another application is footstep planning in robotics~\cite{Deits:2014,Kuindersma:2016}, where $s=S=1$, $x = (\cos(\theta),\sin(\theta))$, and $x$ must satisfy the trigonometric identity $x_1^2 + x_2^2 = 1$.

When $0 < s \leq S$, $\calA$ is a nonconvex set. Moreover, the annulus is not \emph{mixed-integer convex representable}~\cite{Lubin:2017,Lubin:2017a}: that is, there do not exist mixed-integer formulations for the annulus even if we allow the relaxation $R$ to be an arbitrary convex set.

 Foster~\cite{Foster:2013} proposes a disjunctive relaxation for the annulus given as $\hat{\calA} \defeq \bigcup_{i=1}^d P^i$, where each
\begin{equation} \label{eqn:annulus-pieces}
    P^i = \Conv\left(\left\{v^{2i+s-4}\right\}_{s=1}^4\right) \quad \forall i \in \llbracket d \rrbracket
\end{equation}
is a quadrilateral based on the breakpoints
\begin{alignat*}{3}
    v^{2i-1} &= \left(s\cos\left(\frac{2\pi i}{d}\right), s\sin\left(\frac{2\pi i}{d}\right)\right) \quad\quad &\forall i \in \llbracket d \rrbracket& \\
    v^{2i} &= \left(S\sec\left(\frac{2\pi }{d}\right)\cos\left(\frac{2\pi i}{d}\right), S\sec\left(\frac{2\pi }{d}\right)\cos\left(\frac{2\pi i}{d}\right)\right) \quad\quad &\forall i \in \llbracket d \rrbracket&,
\end{alignat*}
where, for notational simplicity, we take $v^{0} \equiv v^{2d}$ and $v^{-1} \equiv v^{2d-1}$. We can in turn represent this disjunctive relaxation through the combinatorial disjunctive constraint given by the family $\TANN_d \defeq (T^i = \{2i+s-4\}_{s=1}^4)_{i=1}^d$. See the right side of Figure~\ref{fig:annulus} for an illustration. For the remainder, when we refer to a formulation of the annulus, we understand this to mean that it is a formulation for the combinatorial disjunctive constraint given by the sets $\calP(\TANN_d)$.

\subsubsection{Small (logarithmic) traditional MIP formulations for the annulus} \label{sec:MIP-annulus}

We start by using Theorem~\ref{thm:general-cdc-characterization} to present new traditional MIP formulations for the annulus. Foster~\cite{Foster:2013} constructs a ``disaggregated logarithmic'' MIP formulation~\cite{Vielma:2010} for $\TANN_d$. This formulation does not take any combinatorial structure of the constraint into account; in our framework, it corresponds to taking each set in $\calT$ as nonintersecting by repeating shared breakpoints (and so $D = \emptyset$). This leads to an increase in the number of components of $\lambda$, as well as a degradation of computational performance relative to logarithmic MIP formulations that use structure~\cite{Huchette:2017,Vielma:2010}.

We start by presenting an ideal logarithmic traditional MIP formulation for the annulus that uses $\lceil\log_2(d)\rceil$ control variables and $2\lceil\log_2(d)\rceil$ general inequality constraints.

\begin{proposition}\label{prop:log-annulus}
    Fix $d = 2^r$ for some $r \in \bbN$. Take the binary reflected Gray encoding $\HRB_d = (h^i)_{i=1}^d \subseteq \{0,1\}^r$, along with $h^0 \equiv h^d$ for notational convenience. Then $(\lambda,z) \in Q(\calP(\TANN_d),\HRB_d)$ if and only if
    \begin{subequations} \label{eqn:annulus-log-form}
    \begin{gather}
        \sum_{i=1}^d \min\{h^{i-1}_k,h^i_k\} (\lambda_{2i-1} + \lambda_{2i}) \leq z_k \quad \forall k \in \llbracket r \rrbracket \\
        \sum_{i=1}^d \max\{h^{i-1}_k,h^i_k\} (\lambda_{2i-1} + \lambda_{2i}) \geq z_k \quad \forall k \in \llbracket r \rrbracket \\
        (\lambda,z) \in \Delta^{2d} \times \bbR^r.
    \end{gather}
    \end{subequations}
\end{proposition}
\proof{}
    The result follows from Theorem~\ref{thm:general-cdc-characterization} after observing that $D = \{i,i+1\}_{i=1}^{d-1} \cup \{1,d\}$ and therefore that $C = \{\pm{\bf e}^k\}_{k=1}^r$, as the binary reflected Gray code is cyclic ($h^{d}-h^1 = {\bf e}^1$).
\qed\endproof

We can also apply Theorem~\ref{thm:general-cdc-characterization} using the zig-zag encoding of Huchette and Vielma~\cite{Huchette:2017} to produce another traditional MIP formulation for the annulus with $\lceil \log_2(d) \rceil$ control variables and $\scrO(\log^2(d))$ general inequality constraints.

\begin{proposition} \label{prop:zig-zag-annulus}
    Fix $d = 2^r$ for some $r \in \bbN$. Take the zig-zag encoding $\HZZ_d = (h^i)_{i=1}^d \subseteq \{0,1\}^r$, along with $h^0 \equiv h^d$ for notational convenience. Then $(\lambda,z) \in Q(\calP(\TANN_d),\HZZ_d)$ if and only if
    \begin{subequations} \label{eqn:annulus-zigzag-form}
    \begin{gather}
        \sum_{i=1}^d \min\{h^{i-1}_k,h^i_k\} (\lambda_{2i-1} + \lambda_{2i}) \leq z_k \quad \forall k \in \llbracket r \rrbracket \\
        \sum_{i=1}^d \max\{h^{i-1}_k,h^i_k\} (\lambda_{2i-1} + \lambda_{2i}) \geq z_k \quad \forall k \in \llbracket r \rrbracket \\
        \sum_{i=1}^d \min\left\{\frac{h^{i-1}_k}{2^{\ell}}-\frac{h^{i-1}_\ell}{2^{k}},\frac{h^{i}_k}{2^{\ell}}-\frac{h^{i}_\ell}{2^{k}}\right\} (\lambda_{2i-1} + \lambda_{2i}) \leq \frac{z_k}{2^\ell} - \frac{z_\ell}{2^k} \quad \forall \{k,\ell\} \in [r]^2 \\
        \sum_{i=1}^d \max\left\{\frac{h^{i-1}_k}{2^{\ell}}-\frac{h^{i-1}_\ell}{2^{k}},\frac{h^{i}_k}{2^{\ell}}-\frac{h^{i}_\ell}{2^{k}}\right\} (\lambda_{2i-1} + \lambda_{2i}) \geq \frac{z_k}{2^\ell} - \frac{z_\ell}{2^k} \quad \forall \{k,\ell\} \in [r]^2 \\
        (\lambda,z) \in \Delta^{2d} \times \bbR^r.
    \end{gather}
    \end{subequations}
\end{proposition}
\proof{}
    The result follows from Theorem~\ref{thm:general-cdc-characterization}. As $D = \{i,i+1\}_{i=1}^{d-1} \cup \{1,d\}$, it follows that $C = \{{\bf e}^k\}_{k=1}^r \cup \{w \equiv (2^{r-1},2^{r-2},\ldots,2^0)\}$. We have that $B = \{{\bf e}^k\}_{k=1}^r$ induce all hyperplanes spanned by the vectors $C \backslash \{w\} = \{{\bf e}^k\}_{k=1}^r$. Now consider each hyperplane spanned by $\hat{C} = \{{\bf e}^k\}_{k \in I} \cup \{w\} \subset C$, where $\hat{I} \subseteq I$. As $|C| = r+1$ and $\dim(C) = r$, we must have $|I| = r-2$, i.e. that there are distinct indices $k,\ell \in \llbracket r \rrbracket \backslash I$ where $I \cup \{k,\ell\} = \llbracket r \rrbracket$. We may then compute that the corresponding hyperplane is given by the normal direction $b^{k,\ell} \defeq 2^{-\ell}{\bf e}^{k} - 2^{-k}{\bf e}^{\ell}$. Therefore, we have that the set $B = \{{\bf e}^k\}_{k=1}^r \cup \{b^{k,\ell}\}_{\{k,\ell\} \in [r]^2}$ suffices for the conditions of Theorem~\ref{thm:general-cdc-characterization}, giving the result.
\qed\endproof
The analysis of Huchette and Vielma~\cite{Huchette:2017} shows that the the zig-zag encoding enjoys the ``incremental branching''' behavior for univariate piecewise linear functions, which leads to computational performance gains relative to the existing logarithmic formulation of Vielma et al.~\cite{Vielma:2010,Vielma:2009a}. Therefore, it may be the case that the zig-zag formulation for the annulus similarly outperforms the logarithmic formulation \eqref{eqn:annulus-log-form}, despite the modest increase in general inequality constraints.

\subsubsection{A very small mixed-integer branching formulation for the annulus}
In addition to the two logarithmic traditional MIP formulations, we can also present a very small mixed-integer branching formulations for the annulus that requires only a constant number of control variables and general inequality constraints.

\begin{proposition} \label{prop:exotic-annulus}
    Take $\HEX_d=(h^i)_{i=1}^d$ as given in \eqref{eqn:exotic-encoding}, along with $h^0 \equiv h^d$ for notational convenience. Then $(\lambda,z) \in Q(\calP(\TANN_d),\HEX_d)$ if and only if
    \begin{subequations} \label{eqn:annulus-constant}
    \begin{align}
        \sum_{i=1}^{d} \min\{h^{i-1}_k,h^{i}_k\}(\lambda_{2i-1}+\lambda_{2i}) &\leq z_k \quad\quad \forall k \in \llbracket 2 \rrbracket \\
        \sum_{i=1}^{d} \max\{h^{i-1}_k,h^{i}_k\}(\lambda_{2i-1}+\lambda_{2i}) &\geq z_k \quad\quad \forall k \in \llbracket 2 \rrbracket \\
        \sum_{i=1}^{d} \min\{w \cdot h^{i-1},w \cdot h^{i}\}(\lambda_{2i-1}+\lambda_{2i}) &\leq w \cdot z \\
        \sum_{i=1}^{d} \max\{w \cdot h^{i-1},w \cdot h^{i}\}(\lambda_{2i-1}+\lambda_{2i}) &\geq w \cdot z \\
        (\lambda,z) \in \Delta^{2d} \times \bbR^{2},
    \end{align}
    where $w = (h^d_2-h^1_2,h^1_1-h^d_1)$.
    \end{subequations}
\end{proposition}
\proof{}
    As $D = \{i,i+1\}_{i=1}^{d-1} \cup \{1,d\}$, then $C \subseteq \{\pm {\bf e}^1, \pm {\bf e}^2, h^d-h^1\}$, and the result immediately follows from Theorem~\ref{thm:general-cdc-characterization}.
\qed\endproof


\subsection{A big-$M$ mixed-integer branching formulation for any disjunctive set}

Our discussion to this point has been restricted to combinatorial disjunctive constraints, for which Theorem~\ref{thm:general-cdc-characterization} gives an explicit geometric construction for ideal formulations. Although combinatorial disjunctive constraints can be used to formulate any disjunctive constraint, this is not always prudent (for example, if the number of extreme points is large). Therefore, we close by presenting a a simple big-$M$ formulation for a disjunctive constraint given by an inequality description.

Consider a generic disjunctive set, where we have an explicit linear inequality description $\calP = (P^i = \{x \in \bbR^n : A^ix \leq b^i\})_{i=1}^d$ for each alternative. In general, it will be very difficult to compute $Q(\calP,H)$ to produce an ideal formulation for the disjunctive constraint. However, we can still apply the standard big-$M$ technique to produce a \emph{non-ideal} mixed-integer branching formulation using only two control variables and a modest number of constraints.

\begin{proposition} \label{prop:big-M}
    Take the family of bounded polyhedra $\calP = (P^i = \{x \in \bbR^n : A^ix \leq b^i \})_{i=1}^d$, where $A^i \in \bbR^{m_i \times n}$ and $b^i \in \bbR^{m_i}$. Take $M^i \in \bbR^{m_i}$ for each $i \in \llbracket d \rrbracket$ such that $M^i_s \geq \max_{x \in \bigcup_{k \neq i} P^k} A^k_s x$ for each $s \in \llbracket m_i \rrbracket$. Then $(x,z) \in \Em(\calP,\HMC_d)$ if and only if
    \begin{subequations}\label{eqn:big-M-formulation}
    \begin{alignat}{2}
        A^i x &\leq b^i + (M^i-b^i)\left( i^2 - 2iz_1 + z_2 \right) \quad &\forall i \in \llbracket d \rrbracket \label{eqn:big-M-formulation-1} \\
        z &\in \HMC_d. \label{eqn:big-M-formulation-2}
    \end{alignat}
    Additionally, we can construct a LP relaxation for a corresponding formulation of $\Em(\calP,\HMC_d)$ by replacing \eqref{eqn:big-M-formulation-2} with the constraint $z \in \Psi_d(1,d)$.
    \end{subequations}
\end{proposition}
\proof{}
    Consider the constraints \eqref{eqn:big-M-formulation-1}, given $z = (i,i^2) \in \HMC_d$. The $j$-th set of constraints in \eqref{eqn:big-M-formulation-1} simplifies to
    \[
        A^jx \leq \begin{cases} b^j + (M^j-b^j)(i^2-2i^2+i^2) = b^i & j = i \\
                                b^j + (M^j-b^j)(j^2-2j \cdot i+i^2) \equiv \alpha^j & \text{o.w.} \end{cases}
    \]
    As $j^2-2i\cdot j+i^2 = (i-j)^2 \geq 1$ for each $i,j \in \bbZ$ with $i \neq j$, we have that $\alpha \geq M^j$. Therefore, given $z=(i,i^2) \in \HMC_d$, $x$ satisfies these constraints if and only if $x \in P^i$.
\qed \endproof

We compare this formulation against a big-$M$ traditional MIP formulation~\cite{Vielma:2015}. Both require $\sum_{i=1}^d m_i$ general inequality constraints, along with $\scrO(d)$ additional constraints to describe either $\Psi_d(1,d)$ or variable bounds on binary variables. However, formulation \eqref{eqn:big-M-formulation} requires only two control variables, compared to the $\lceil \log_2(d) \rceil$ binary control variables needed for a traditional big-$M$ MIP formulation.





\section{Omitted Proofs}

\subsection{Proof of Proposition~\ref{easycombinatorial}} \label{app:prove-cdc-properties}

\proof{}
    
    The ``if'' direction will follow if we can show that $Q(\calP(\calT),H) \cap (\bbR^n \times H) = \Em(\calP(\calT),H)$, as then $Q(\calP(\calT,H))$ is the relaxation for a valid formulation for $\bigcup_{i=1}^d P(T^i)$. To see this, consider each $h^i \in H$, along with the associated $\Slice(h^i) = \{x : (x,h^i) \in Q(\calP(\calT),H)\}$. Clearly $\Slice(h^i) \supseteq P(T^i)$ from the definition of $\Em(\calP(\calT),H)$, and so the result will follow if we show that $\Slice(h^i) \subseteq P(T^i)$. 
    
    Take some $\hat{x} \in \Slice(h^i)$; then $(\hat{x},h^i) \in Q(\calP(\calT),H)$ necessarily. Therefore, it is possible to express $(\hat{x},h^i)$ as a convex combination of points in $\Em(\calP(\calT),H)$. Equivalently, there exists some $\lambda \in \Delta^d$ and some points $\tilde{x}^j \in P(T^j)$ for each $j \in \llbracket d \rrbracket$ such that $(\hat{x},h^i) = \sum_{j=1}^d \lambda_j (\tilde{x}^j,h^j)$. As the codes in $H$ are in convex position, then $h^i = \sum_{j=1}^d \lambda_j h^j$ implies that $\lambda = {\bf e}^i$, giving the result.
    
    For the ``only if'' direction, we start by showing that any formulation for $\bigcup_{i=1}^d P(T^i)$ with relaxation $R$ must necessarily satisfy the property that for each $i \in \llbracket d \rrbracket$, $\Slice(h^i) = \{x : (x,h^i) \in R\} = P(T^{j_i})$ for some index $j_i \in \llbracket d \rrbracket$. For each such $i$, clearly $\Slice(h^i) \subseteq P(T^{j_i})$ for some $j_i \in \llbracket d \rrbracket$, since otherwise the formulation cannot be valid. By our nonredundancy assumption and the pigeon-hole principle, we can assign each $i \in \llbracket d \rrbracket$ a unique index $j_i \in \llbracket d \rrbracket$; w.l.o.g., take $j_i = i$ for each $i \in \llbracket d \rrbracket$.
    
    Now presume for contradiction that $H$ is not in convex position. This presumption, along with the fact that $H$ is a set of distinct vectors, implies there exists a code (w.l.o.g. $h^1$) where $h^1 = \sum_{i=1}^d \lambda_i h^i$ for some $\lambda \in \Delta^d$ with $\lambda_1 = 0$ and at least two fractional components: w.l.o.g., $0 < \lambda_2, \lambda_3 < 1$.

    From the validity of the formulation, it necessarily follows that
    \begin{align*}
        \sum\nolimits_{j=1}^d \lambda_j (P(T^j) \times \{h^j\}) &= \sum\nolimits_{j=1}^d (\lambda_j P(T^j) \times \{h^1\}) \\
        &\subseteq R \cap (\bbR^{n} \times \{h^1\}) \\
        &= P(T^1) \times \{h^1\},
    \end{align*}
    where the inclusion follows from the convexity of $R$. Therefore, $P(T^1) \supseteq \sum_{j=2}^d \lambda_j P(T^j)$. Consider some $v \in T^2$. Since ${\bf e}^v \in P(T^2)$ and each point on the unit simplex is nonnegative, we must have that $\tilde{\lambda} \in P(T^1)$ for some $\tilde{\lambda} \in \Delta^n$ with $ \tilde{\lambda}_v \geq \lambda_v > 0$. Therefore, we must have ${\bf e}^v \in P(T^1)$, and hence $v \in T^1$. Repeating this for each $v \in T^2$, we conclude that $T^1 \supseteq T^2$. However, this contradicts the irredundancy assumption, and so we must have that $H$ is in convex position.
    
    The lower bound on the number of control variables for hole-free encodings follows as $\Conv(H)$ has at most $2^r$ extreme points \cite[Proposition 3]{Celaya:2015}, implying that $H$ has at most $d=2^r$ elements.
\qed \endproof

\subsection{Proof of Theorem~\ref{thm:general-cdc-characterization}} \label{sec:prove-cdc-characterization}

\proof{}
    It is straightforward to show the ``only if'' direction; indeed, the essential argument already presented in the proof of Proposition~\ref{prop:general-cdc-with-general-2D-embedding} suffices. Therefore, we focus on the ``if'' direction.

    From the connectivity assumption on $D$, it is possible to show that $\calL = \aff(H) - h^1$, where the choice of $h^1$ to subtract from $\aff(H)$ was arbitrary. Now, let $F$ be a facet of $Q(\calP,H)$. By possibly adding or subtracting multiples of $\sum_{i=1}^n \lambda_i = 1$ and the equations defining $\aff(H)$, we may assume w.l.o.g. that $F$ is induced by $\tilde{a} \cdot \lambda \leq \tilde{b} \cdot y$ for some $(\tilde{a},\tilde{b}) \in \bbR^{n+r}$.  If $B = \ext(Q(\calP,H))$ is the set of all extreme points, we will have that $F$ is supported by some strict nonempty subset $\tilde{B} \subsetneq B$. Take $\tilde{D} = \{\{i,j\} \in D : \exists v \in \llbracket n \rrbracket \text{ s.t. } ({\bf e}^v,h^i),({\bf e}^v,h^j) \in \tilde{B}\}$ and $\tilde{C} = \{c^{i,j} \in C : \{i,j\} \in \tilde{D}\}$. In particular, we see that $\tilde{b} \cdot c^{i,j} = 0$ for each $c^{i,j} \in \tilde{C}$, as if $\{i,j\} \in \tilde{D}$, this implies that there is some $v \in \llbracket n \rrbracket$ whereby $\tilde{a} \cdot {\bf e}^v = \tilde{b} \cdot h^i = \tilde{b} \cdot h^j$.

    \paragraph{\underline{Case 1: $\dim(\tilde{C}) = \dim(C)$}}
    In this case, we show that $F$ corresponds to a variable bound on a single component of $\lambda$. As $\tilde{C} \subseteq C$ and $\dim(\tilde{C}) = \dim(C)$, we conclude that $\Span(\tilde{C}) = \Span(C) \equiv \calL$. Then $\tilde{b} \in \calL^\perp$, as $\tilde{b} \perp \tilde{C}$. Furthermore, $\calL$ is the linear space parallel to $\aff(H)$. Therefore, we can w.l.o.g. presume that $\tilde{b} = {\bf 0}^r$, as \eqref{eqn:general-V-formulation-2} constrains $z \in \aff(H)$.

    We observe that $\tilde{a} \neq {\bf 0}^n$, as otherwise this would correspond to the vacuous inequality $0 \leq 0$, which is not a proper face. We now show that $\tilde{a}$ has exactly one nonzero element. Assume for contradiction that this is not the case, and w.l.o.g. $\tilde{a}_1, \tilde{a}_2 < 0$ (any strictly positive components will not yield a valid inequality for $B$). This would imply that $({\bf e}^1,h^j) \in B^\star$ for each $j \in \llbracket d \rrbracket$ such that $1 \in T^j$, and similarly that $({\bf e}^2,h^j) \in B^\star$ for each $j \in \llbracket d \rrbracket$ wherein $2 \in T^j$. However, we could then perform the simple tilting $\tilde{a}_2 \leftarrow 0$ to construct a distinct face with strictly larger support, as now $({\bf e}^2,h^j)$ is supported by the corresponding face for each $j$ such that $2 \in T^j$. Furthermore, as this new constraint does not support $({\bf e}^1,h^j)$ for each $j$ such that $1 \in T^j$, the new face is proper, and thus contradicts the original face $F$ being a facet. Therfore, we can normalize the coefficients to $\tilde{a} = -{\bf e}^1$, giving a variable bound constraint on a component of $\lambda$ which appears in the restriction $\lambda \in \Delta^n$ in \eqref{eqn:general-V-formulation-2}.

    \paragraph{\underline{Case 2: $\dim(\tilde{C}) = \dim(C) - 1$}}
    The fact that $b \perp \tilde{C}$, along with the dimensionality of $\tilde{C}$, implies that $M(\tilde{b};\calL) = \Span(\tilde{C})$ is a hyperplane in $\calL$. This means we can assume w.l.o.g. that $\tilde{b} = s b^i$ for some $i \in \llbracket \Gamma \rrbracket$ and $s \in \{-1,+1\}$. We then compute for each $v \in \llbracket n \rrbracket$ that either $a_v = \min_{j : v \in T^j}\{b^i \cdot h^j\}$ if $s=+1$, or $a_v = -\max_{j : v \in T^j}\{b^i \cdot h^j\}$ if $s = -1$.

    \paragraph{\underline{Case 3: $\dim(\tilde{C}) < \dim(C)-1$}}
    We will show that this case cannot occur if $F$ is a general inequality facet. In fact, observe that if, w.l.o.g. ${\bf e}^1 \not\in \Proj_\lambda(\tilde{B})$, then $\tilde{a} \cdot \lambda \leq \tilde{b} \cdot z$ is either equal to, or dominated by, the variable bound $\lambda_1 \geq 0$. Therefore, we assume that $\Proj_\lambda(\tilde{B}) = \{{\bf e}^i\}_{i=1}^n$ for the remainder.

    Presume for contradiction that it is indeed the case that $F$ is a facet and $\dim(\tilde{C}) < \dim(C)-1$. As $F$ is a proper face, we know that there is some point in $B$ not supporting $F$, w.l.o.g. $({\bf e}^1,h^1) \in B \backslash \tilde{B}$. We will take all the remaining extreme points as $B^\star = B \backslash (\tilde{B} \cup \{({\bf e}^1,h^1)\})$.

    First, we show that $B^\star \neq \emptyset$. If this where not the case, then $\tilde{B} = B \backslash \{({\bf e}^1,h^1)\}$ necessarily, and this implies that $i=1$ for each $\{i,j\} \in D \backslash \tilde{D}$ (recall that $i < j$ notationally). Furthermore, $T^1 \cap T^j \subseteq \{1\}$ for each $j \in \llbracket 2,d \rrbracket$, else $c^{1,v} \in \tilde{C}$ and $\{1,v\} \in \tilde{D}$. Therefore, as $D$ is connected by assumption, $\tilde{D}$ is ``nearly connected'' in the sense that $G = (\llbracket 2,d \rrbracket, \{\{i,j\} \in \tilde{D} : i \neq 1, j \neq 1 \})$ is a connected graph. By the same argument as in the beginning of the proof, we conclude that $\Span(\tilde{C}) \supseteq \aff(\{h^i\}_{i=2}^d) - h^2$. However, this would imply that $\dim(\tilde{C}) \geq \dim(\aff(\{h^i\}_{i=2}^d) - h^2) \geq \dim(\aff(\{h^i\}_{i=1}^d) - h^2) - 1 = \dim(\calL) - 1 = \dim(C)-1$, a contradiction of our dimensionality assumption. Therefore, we conclude that $B^\star \neq \emptyset$.

    We now define the cone
    \[
        K = \left\{(a,b) \in \bbR^n \times \calL : a \cdot {\bf e}^v \leq b \cdot h^j \: \forall ({\bf e}^v,h^j) \in B^\star \right\}
    \]
    and the linear space
    \[
        L = \left\{ (a,b) \in \bbR^n \times \calL : a \cdot {\bf e}^v = b \cdot h^j \: \forall ({\bf e}^v,h^j) \in \tilde{B} \right\}.
    \]
    Furthermore, we see that the inequalities defining $K$ cannot be implied equalities. Therefore, as $(\tilde{a},\tilde{b}) \in K$ and this point strictly satisfies each inequality of $K$ indexed by $B^\star$, we conclude that $K$ is full-dimensional in $\bbR^n \times \calL$, and that $(\tilde{a},\tilde{b}) \in \relint(K)$.

    Next, we show that $\dim(L) > 1$. To show this, we start by instead studying $L' = \{ b \in \calL : b \cdot c = 0 \: \forall c \in \tilde{C} \}$. We can readily observe that $L' = \Proj_b(L)$. Furthermore, as $\Proj_a(\tilde{B}) = \{{\bf e}^i\}_{i=1}^n$ from the argument at the beginning of the case, we conclude that the set $\{a : (a,b) \in L\}$ is a singleton. In other words, the values for $a$ are completely determined by the values for $b$ in $L$. From this, we conclude that $\dim(L) = \dim(L')$. From the definition of $L'$, we see that $L'$ and $\Span(\tilde{C})$ form an orthogonal decomposition of $\calL$. Therefore, $\dim(\calL) = \dim(L') + \dim(\tilde{C})$. Recalling that $\dim(\calL) = \dim(C)$, and that we are assuming that $\dim(\tilde{C}) < \dim(C)-1$, we have that $\dim(L) = \dim(L') = \dim(\calL) - \dim(\tilde{C}) = \dim(C) - \dim(\tilde{C}) > 1$, giving the result.

    We now show that $K \cap L$ is pointed. To see this, presume for contradiction that there exists a nonzero $(\hat{a},\hat{b})$ such that $(\hat{a},\hat{b}),(-\hat{a},-\hat{b}) \in K\cap L$. However, this would imply that $\hat{a} \cdot {\bf e}^v = \hat{b} \cdot h^j$ for all $({\bf e}^v,h^j) \in \tilde{B} \cup B^\star$. Because $B^\star\neq\emptyset$, this implies that $\hat{a} \cdot \lambda \leq \hat{b} \cdot z$ is a face strictly containing the facet $F$, and so must be a non-proper face (i.e. it is additionally supported by $({\bf e}^1,h^1)$ and hence by every point in $B$). However, this would imply that $\hat{b} \cdot c = 0$ for all $c \in C$, and as $\calL = \Span(C)$, this would necessitate that $\hat{b} \in \calL^\perp$. As $\hat{b} \in \calL$ from the definition of $K$, it follows that $\hat{b} = {\bf 0}^r$. However, this would imply that $\hat{a} \cdot \lambda = 0$ is valid for $B$, which cannot be the case unless $\hat{a} = {\bf 0}^n$, a contradiction. Therefore, $K \cap L$ is pointed.

    As $\dim(L) > 1$, we can take some two-dimensional linear subspace $L^2 \subseteq L$ such that $(\tilde{a},\tilde{b}) \in L^2$. As $(\tilde{a},\tilde{b}) \in L \cap \relint(K)$, it follows that $(\tilde{a},\tilde{b}) \in L^2 \cap \relint(K)$ as well. Similarly, as $K \cap L$ is pointed, it follows that $K^2 = L^2 \cap K$ is pointed as well. Furthermore, as $K$ is full-dimensional in $\bbR^n \times \calL$, $K^2$ is full-dimensional in $L^2 \subset \bbR^n \times \calL$ (i.e. 2-dimensional).Therefore, a minimal description for it includes the equalities that define $L^2$, along with exactly two nonempty-face-inducing inequality constraints from the definition of $K$. Add the single strict inequality $\hat{K}^2 = K^2 \cap \{(a,b) \in \bbR^n \times \calL : a \cdot {\bf e}^1 < b \cdot h^1\}$. As $\tilde{a} \cdot {\bf e}^1 < \tilde{b} \cdot h^1$ and $(\tilde{a},\tilde{b}) \in K^2$, it follows that $\hat{K}^2$ is nonempty and also 2-dimensional, and can be described using only the linear equations defining $L^2$, the strict inequality $a \cdot {\bf e}^1 < b \cdot h^1$, and at least one (and potentially two) of the inequalities previously used to describe $K^2$. Select one of the defining nonempty-face-inducing inequalities given by $a \cdot {\bf e}^v \leq b \cdot h^j$, where $({\bf e}^v,h^j) \in B^\star$.

    Now construct the restriction $S = \{(a,b) \in \hat{K}^2 : a \cdot {\bf e}^v = b \cdot h^j\}$. As $a \cdot {\bf e}^v \leq b \cdot h^j$ induces a non-empty face on the cone $\hat{K}^2$, $S$ is nonempty. Furthermore, we see that any $(\hat{a},\hat{b}) \in S$ will correspond to a valid inequality $\hat{a} \cdot \lambda \leq \hat{b} \cdot z$ for $B$ with strictly greater support than our original face $\tilde{a} \cdot \lambda \leq \tilde{b} \cdot z$. In particular, we see that $({\bf e}^v,h^j) \in B^\star$, i.e. $\tilde{a} \cdot {\bf e}^v < \tilde{b} \cdot h^j$, but by construction $\hat{a} \cdot {\bf e}^v = \hat{b} \cdot h^j$. Additionally, since $\hat{a} \cdot {\bf e}^1 < \hat{b} \cdot h^1$, the corresponding face is proper, which implies that $F$ cannot be a facet, a contradiction.\qed
\endproof

\begin{acknowledgements}
    This material is based upon work supported by the National Science Foundation under Grant CMMI-1351619. We thank Matthias K\"oppe for suggesting a potential connection between the embedding object and points along the moment curve.
\end{acknowledgements}

\bibliographystyle{spmpsci}
\bibliography{master.bib}

\end{document}